\documentclass[11pt]{article}
\usepackage{fullpage}
\usepackage[utf8]{inputenc}
\usepackage[T1]{fontenc}
\usepackage{graphicx}
\usepackage{subcaption}

\usepackage{floatrow}
\newfloatcommand{capbtabbox}{table}[][\FBwidth]
\usepackage{textcomp}
\usepackage{amsbsy}
\usepackage{latexsym}
\usepackage[mathscr]{eucal}
\usepackage{amsfonts}
\usepackage{amssymb}
\usepackage{mathtools}
\usepackage[usenames]{color}
\usepackage[dvipsnames]{xcolor}

\usepackage{algorithm, algcompatible}
\renewcommand{\COMMENT}[2][.6\linewidth]{%
  \leavevmode\hfill\makebox[#1][l]{//~#2}}
\algnewcommand\algorithmicto{\textbf{to}}
\algnewcommand\RETURN{\State \textbf{return} }

\graphicspath{{img/}}

\def\wnew{\color{black}}

\begin{document}
\bibliographystyle{plain}
\title
{
Optimal reduced model algorithms for data-based state estimation}
\author{ 
Albert Cohen, Wolfgang Dahmen, Ron DeVore, Jalal Fadili, Olga Mula and James Nichols
\thanks{%
This research was supported by the
Institut Universitaire de France;
the ERC Adv grant BREAD;
the EMERGENCES grant of the Paris City Council ``Models and Measures''  and
the NSF grants   DMS 15-21067,   DMS  18-17603, ONR grants
  N00014-17-1-2908, N00014-16-1-2706 (RD). A portion of this research was completed while A.C. W.D.,  R.D. (Simon Fellow),   O.M. and J.N. were supported as  visitors of the Isaac Newton Institute at Cambridge University.
}  }

\hbadness=10000
\vbadness=10000
\newtheorem{lemma}{Lemma}[section]
\newtheorem{prop}[lemma]{Proposition}
\newtheorem{cor}[lemma]{Corollary}
\newtheorem{theorem}[lemma]{Theorem}
\newtheorem{remark}[lemma]{Remark}
\newtheorem{example}[lemma]{Example}
\newtheorem{definition}[lemma]{Definition}
\newtheorem{proper}[lemma]{Properties}
\newtheorem{assumption}[lemma]{Assumption}
%
\newenvironment{disarray}{\everymath{\displaystyle\everymath{}}\array}{\endarray}

\def\vp{\varphi}
\def\<{\langle}
\def\>{\rangle}
\def\t{\tilde}
\def\i{\infty}
\def\e{\varepsilon}
\def\sm{\setminus}
\def\nl{\newline}
\def\o{\overline}
\def\wt{\widetilde}
\def\wh{\widehat}
\def\cT{{\cal T}}
\def\cA{{\cal A}}
\def\cI{{\cal I}}
\def\cV{{\cal V}}
\def\cB{{\cal B}}
\def\cF{{\cal F}}

\def\cR{{\cal R}}
\def\cD{{\cal D}}
\def\cP{{\cal P}}
\def\cJ{{\cal J}}
\def\cM{{\cal M}}
\def\cO{{\cal O}}
\def\Chi{\raise .3ex
\hbox{\large $\chi$}} \def\vp{\varphi}
\def\lsima{\hbox{\kern -.6em\raisebox{-1ex}{$~\stackrel{\textstyle<}{\sim}~$}}\kern -.4em}
\def\lsim{\hbox{\kern -.2em\raisebox{-1ex}{$~\stackrel{\textstyle<}{\sim}~$}}\kern -.2em}
\def\[{\Bigl [}
\def\]{\Bigr ]}
\def\({\Bigl (}
\def\){\Bigr )}
\def\[{\Bigl [}
\def\]{\Bigr ]}
\def\({\Bigl (}
\def\){\Bigr )}
\def\L{\pounds}
\def\pr{{\rm Prob}}
\newcommand{\cs}[1]{{\color{magenta}{#1}}}
\def\ds{\displaystyle}
\def\ev#1{\vec{#1}}     
\newcommand{\lt}{\ell^{2}(\nabla)}
\def\Supp#1{{\rm supp\,}{#1}}
\def\R{\mathbb{R}}
\def\E{\mathbb{E}}
\def\nl{\newline}
\def\T{{\relax\ifmmode I\!\!\hspace{-1pt}T\else$I\!\!\hspace{-1pt}T$\fi}}
\def\N{\mathbb{N}}
\def\Z{\mathbb{Z}}
\def\N{\mathbb{N}}
\def\Zd{\Z^d}
\def\Q{\mathbb{Q}}
\def\C{\mathbb{C}}
\def\Rd{\R^d}
\def\gsim{\mathrel{\raisebox{-4pt}{$\stackrel{\textstyle>}{\sim}$}}}
\def\sime{\raisebox{0ex}{$~\stackrel{\textstyle\sim}{=}~$}}
\def\lsim{\raisebox{-1ex}{$~\stackrel{\textstyle<}{\sim}~$}}
\def\div{\mbox{ div }}
\def\M{M}  \def\NN{N}                  
\def\L{{\ell}}               
\def\Le{{\ell^1}}            
\def\Lz{{\ell^2}}
\def\Let{{\tilde\ell^1}}     
\def\Lzt{{\tilde\ell^2}}
\def\Ltw{\ell^\tau^w(\nabla)}
\def\t#1{\tilde{#1}}
\def\la{\lambda}
\def\La{\Lambda}
\def\ga{\gamma}
\def\BV{{\rm BV}}
\def\Ga{\eta}
\def\al{\alpha}
\def\cZ{{\cal Z}}
\def\cA{{\cal A}}
\def\cU{{\cal U}}
\def\ms{{\rm ms}}
\def\wc{{\rm wc}}
\def\argmin{\mathop{\rm argmin}}
\def\argmax{\mathop{\rm argmax}}
\def\prob{\mathop{\rm prob}}
\def\A{\mathop{\rm Alg}}
\def\cen{{\rm cen}}

\def \bphi{{\bf\phi}}

\def\cO{{\cal O}}
\def\cA{{\cal A}}
\def\cC{{\cal C}}
\def\cS{{\cal F}}
\def\bu{{\bf u}}
\def\bz{{\bf z}}
\def\bZ{{\bf Z}}
\def\bI{{\bf I}}
\def\cE{{\cal E}}
\def\cD{{\cal D}}
\def\cG{{\cal G}}
\def\cI{{\cal I}}
\def\cJ{{\cal J}}
\def\cM{{\cal M}}
\def\cN{{\cal N}}
\def\cT{{\cal T}}
\def\cU{{\cal U}}
\def\cV{{\cal V}}
\def\cW{{\cal W}}
\def\cL{{\cal L}}
\def\cB{{\cal B}}
\def\cG{{\cal G}}
\def\cK{{\cal K}}
\def\cS{{\cal S}}
\def\cP{{\cal P}}
\def\cQ{{\cal Q}}
\def\cR{{\cal R}}
\def\cU{{\cal U}}
\def\bL{{\bf L}}
\def\bl{{\bf l}}
\def\bK{{\bf K}}
\def\bQ{{\bf Q}}
\def\bC{{\bf C}}
\def\X{X\in\{L,R\}}
\def\ph{{\varphi}}
\def\D{{\Delta}}
\def\H{{\cal H}}
\def\bM{{\bf M}}
\def\bx{{\bf x}}
\def\bj{{\bf j}}
\def\bG{{\bf G}}
\def\bS{{\bf S}}
\def\bP{{\bf P}}
\def\bW{{\bf W}}
\def\bX{{\bf X}}
\def\bT{{\bf T}}
\def\bV{{\bf V}}
\def\bv{{\bf v}}
\def\bt{{\bf t}}
\def\bz{{\bf z}}
\def\bw{{\bf w}}
\def \mvn {{\rm mvn}}
\def\rhom{{\rho^m}}
\def\diff{\hbox{\tiny $\Delta$}}
\def\EE{{\rm Exp}}
\def\lll{\langle}
\def\argmin{\mathop{\rm argmin}}
\def\argmax{\mathop{\rm argmax}}
\def\dJ{\nabla}
\newcommand{\ba}{{\bf a}}
\newcommand{\bb}{{\bf b}}
\newcommand{\bc}{{\bf c}}
\newcommand{\bd}{{\bf d}}
\newcommand{\bs}{{\bf s}}
\newcommand{\bff}{{\bf f}}
\newcommand{\bp}{{\bf p}}
\newcommand{\bg}{{\bf g}}
\newcommand{\by}{{\bf y}}
\newcommand{\br}{{\bf r}}
\newcommand{\be}{\begin{equation}}
\newcommand{\ee}{\end{equation}}
\newcommand{\bea}{$$ \begin{disarray}{lll}}
\newcommand{\eea}{\end{disarray} $$}
\def \Vol{\mathop{\rm  Vol}}
\def \mes{\mathop{\rm mes}}
\def\rad{\mathop{\rm rad}}
\def \Prob{\mathop{\rm  Prob}}
\def \exp{\mathop{\rm    exp}}
\def \sign{\mathop{\rm   sign}}
\newcommand{\mult}{\mathop{\rm   mult}}
\newcommand{\one}{\mathop{\rm   one}}
\def \wca{{\rm    wca}}
\def \msa{{\rm    msa}}
\def \sp{\mathop{\rm   span}}
\def \vphi{{\varphi}}
\def \csp{\overline \mathop{\rm   span}}
%
%
\newcommand{\beqn}{\begin{equation}}
\newcommand{\eeqn}{\end{equation}}
\def\beginproof{\noindent{\bf Proof:}~ }
\def\endproof{\hfill\rule{1.5mm}{1.5mm}\\[2mm]}

\newcommand{\utr}{u^{\rm true}}
\newcommand{\Cor}{\kappa}

\newcommand{\prox}{\mathrm{prox}}
\newcommand{\epi}{\text{epi}}
\newcommand{\pa}[1]{\left(#1\right)}
\newcommand{\xob}{x_{\mathrm{ob}}}
\newcommand{\xsol}{x^{\star}}
\newcommand{\xbar}{\bar{x}}
\newcommand{\xkm}{x_{k-1}}
\newcommand{\xk}{x_{k}}
\newcommand{\xkp}{x_{k+1}}
\newcommand{\xbarkp}{\bar{x}_{k+1}}
\newcommand{\tk}{t_{k}}
\newcommand{\tkp}{t_{k+1}}
\newcommand{\tbarkp}{\bar{t}_{k+1}}
\newcommand{\xtk}{\pa{\xk,\tk}}
\newcommand{\xtkp}{\pa{\xkp,\tkp}}
\newcommand{\xtbarkp}{\pa{\xbarkp,\tbarkp}}

\newcommand{\vsol}{v^{\star}}
\newcommand{\vbar}{\bar{v}}
\newcommand{\vkm}{v_{k-1}}
\newcommand{\vik}{v_{i,k}}
\newcommand{\vikp}{v_{i,k+1}}
\newcommand{\xiik}{\xi_{i,k}}
\newcommand{\xiikp}{\xi_{i,k+1}}
\newcommand{\vxik}{\pa{\vik,\xiik }_{i \in [p]}}
\newcommand{\vxiok}{\pa{\vik,\xiik}}
\newcommand{\vxikp}{\pa{\vikp,\xiikp}_{i \in [p]}}
\newcommand{\vxiokp}{\pa{\vikp,\xiikp}}

\newenvironment{Proof}{\noindent{\bf Proof:}\quad}{\endproof}

\renewcommand{\theequation}{\thesection.\arabic{equation}}
\renewcommand{\thefigure}{\thesection.\arabic{figure}}

\makeatletter
\@addtoreset{equation}{section}
\makeatother

\newcommand\abs[1]{\left|#1\right|}
\newcommand\clos{\mathop{\rm clos}\nolimits}
\newcommand\trunc{\mathop{\rm trunc}\nolimits}
\renewcommand\d{d}
\newcommand\dd{\mathrm d}
\newcommand\diag{\mathop{\rm diag}}
\newcommand\dist{\mathop{\rm dist}}
\newcommand\diam{\mathop{\rm diam}}
\newcommand\cond{\mathop{\rm cond}\nolimits}
\newcommand\eref[1]{{\rm (\ref{#1})}}
\newcommand{\iref}[1]{{\rm (\ref{#1})}}
\newcommand\Hnorm[1]{\norm{#1}_{H^s([0,1])}}
\def\int{\intop\limits}
\renewcommand\labelenumi{(\roman{enumi})}
\newcommand\lnorm[1]{\norm{#1}_{\ell^2(\Z)}}
\newcommand\Lnorm[1]{\norm{#1}_{L_2([0,1])}}
\newcommand\LR{{L_2(\R)}}
\newcommand\LRnorm[1]{\norm{#1}_\LR}
\newcommand\Matrix[2]{\hphantom{#1}_#2#1}
\newcommand\norm[1]{\left\|#1\right\|}
\newcommand\ogauss[1]{\left\lceil#1\right\rceil}
\newcommand{\QED}{\hfill
\raisebox{-2pt}{\rule{5.6pt}{8pt}\rule{4pt}{0pt}}%
  \smallskip\par}
\newcommand\Rscalar[1]{\scalar{#1}_\R}
\newcommand\scalar[1]{\left(#1\right)}
\newcommand\Scalar[1]{\scalar{#1}_{[0,1]}}
\newcommand\Span{\mathop{\rm span}}
\newcommand\supp{\mathop{\rm supp}}
\newcommand\ugauss[1]{\left\lfloor#1\right\rfloor}
\newcommand\with{\, : \,}
\newcommand\Null{{\bf 0}}
\newcommand\bA{{\bf A}}
\newcommand\bB{{\bf B}}
\newcommand\bR{{\bf R}}
\newcommand\bD{{\bf D}}
\newcommand\bE{{\bf E}}
\newcommand\bF{{\bf F}}
\newcommand\bH{{\bf H}}
\newcommand\bU{{\bf U}}
\newcommand\cH{{\cal H}}
\newcommand\sinc{{\rm sinc}}
\def\enorm#1{| \! | \! | #1 | \! | \! |}

\newcommand{\dm}{\frac{d-1}{d}}

\let\bm\bf
\newcommand{\balpha}{{\mbox{\boldmath$\alpha$}}}
\newcommand{\bbeta}{{\mbox{\boldmath$\beta$}}}
\newcommand{\bal}{{\mbox{\boldmath$\alpha$}}}
\newcommand{\bbi}{{\bm i}}

\newcommand{\test}{\text{test}}
\newcommand{\greedy}{\text{greedy}}

\def\nnew{\color{black}}
\def\mnew{\color{black}}

\newcommand{\dI}{\Delta}
\maketitle
\date{}
 
\begin{abstract}
Reduced model spaces, such as reduced basis and polynomial chaos, are linear spaces
$V_n$ of finite dimension $n$ which are designed for the efficient approximation of certain families of parametrized
PDEs in a Hilbert space $V$. The manifold $\cM$ that gathers
the solutions of the PDE for all admissible parameter values 
is globally approximated by the space $V_n$ with some controlled accuracy $\e_n$,
which is typically much smaller than when using standard approximation spaces of the same dimension such as finite elements.
Reduced model spaces have also been proposed in \cite{MPPY} as a vehicle to design a simple 
linear recovery algorithm of the state $u\in \cM$ corresponding to a particular solution instance
when the values of parameters are unknown but a set of data  is given by $m$ linear measurements of the state. 
The measurements are of the form $\ell_j(u)$, $j=1,\dots,m$, where the $\ell_j$ are linear functionals on $V$. 
The analysis of this approach in \cite{BCDDPW2} shows that the recovery error is 
bounded by $\mu_n \e_n$, where $\mu_n=\mu(V_n,W)$ is the inverse of an inf-sup 
constant that describe the angle between  $V_n$ and the space 
$W$ spanned by the Riesz representers of $(\ell_1,\dots,\ell_m)$. 
A reduced model space which is efficient for approximation might thus be 
ineffective for recovery if $\mu_n$ is large or infinite. In this paper, we discuss the
existence and effective construction of an optimal reduced model space for this recovery method.
We extend our search to affine spaces which are better adapted than
linear spaces for various purposes. Our basic observation is that this problem is equivalent 
to the search of an optimal affine algorithm
for the recovery of $\cM$ in the worst case error sense. This allows us to peform our
search by a convex optimization procedure. Numerical tests illustrate that the reduced model
spaces constructed from our approach perform better than the classical reduced basis spaces.
\end{abstract}

\section{Introduction}

{\wnew
\subsection{Background and context}
State estimation refers to the general problem of approximately recovering the 
true state of a physical system of interest from incomplete data. This task is 
ubiquitous in applied sciences and engineering. One can draw a distinction between
two different application scenarios:
\begin{enumerate}
\item
The physical properties of the states, sometimes referred to as background information, 
are approximately modeled by a nonlinear dynamical system which by itself is neither sufficiently accurate nor stable to
warrant reliable predictions. This is typically the case
 in weather prediction, climatology, or generally in atmospheric research.
One therefore utilizes observational data to {\em correct} the model-based predictions, ideally in real time.
Such correction mechanisms are often based on statistical hypotheses such as Gaussianity of error distributions.
One important approach is {\em ensemble Kalman filtering}  which can be viewed as a recursive Bayesian estimation based on 
Monte Carlo approximations to the first and second moments of the error distributions \cite{LSZ}. A second class of methods
are so called {\em variational data assimilation} 
schemes like 3D-VAR or 4D-VAR \cite{Lor}.
The state predicted by the model 
is then corrected by minimizing a quadratic cost 
functional involving inverses covariance matrices for the background model error and observation error.
In this first scenario, the error bounds between the exact and estimated state
are typically expressed in an average sense, based on the accepted 
simplified statistical model assumptions.
\item
The physical states of interest are reliably described 
in terms of a {\em parameter dependent} family of PDEs which for each
parameter instance can be computed within a desired target accuracy. 
The states are therefore elements of the associated {\em solution manifold}
that consists of all solutions to the PDE as parameters vary.
The task is then to  estimate a state in (or near to) the solution manifold 
from only  a finite number of measurements generated through a {\em fixed} number of sensors. 
 A classical example is to estimate a pressure field of a porous media flow from a finite number of 
 pressure head measurements. The parametric model then could arise from
a Karhunen-Lo\`{e}ve expansion of a random field of permeability coefficients in Darcy's law,
and may thus involve a large or even infinite number of parameters. 
Problems of this type have been investigated over the past decade in the context of Uncertainty Quantification. 
Again, a prominent approach is Bayesian inversion where prior information is given in terms of a probability distribution
for the parameter, inducing a probability distribution for the state \cite{St}. The objective is then to approximate
the posterior probability distribution of states given the data. High dimensionality renders such methods computationally expensive.
Alternatively, state estimation can be formulated  as a constrained optimization problem. For instance, one could minimize  
 the deviation of state measurements over the solution manifold asking for probabilistic or deterministic error bounds.
  In practice, one typically chooses first a sufficiently fine discretization of the {\em high fidelity} continuum model  which then
gives rise to a large scale (discrete) constrained non-convex optimization problem that needs to be solved for each instance of  data.
Ill-posedness of the inversion task necessitates adding regularization terms  which  introduce  a further  
ambiguous bias.  Reduced models are used to alleviate
the possibly prohibitive cost of the numerous forward simulations that are needed in the descent method. A central issue is then to judiciously switch between
the high fidelty model, given in terms of the fine scale discretization, and the low fidelity reduced model, see \cite{willcox}.
\end{enumerate}

In this article, we consider scenario (ii) but pursue a different approach  taking up on recent work in \cite{BCDDPW2,MPPY}. 
Although it can be formulated without any reference to a statistical model, it has conceptual similarities 
with the 3D and 4D-Var variational approach invoked for scenario (i), see \cite{KBGV} and \cite{Ta}
for such connexions.  In contrast to Bayesian inversion, this approach yields deterministic 
error bounds expressed in a worst case sense over the solution manifold, which is the 
primary interest in this paper.  
Specifically, we follow   
\cite{BCDDPW2} and formulate state estimation  as an  {\em optimal recovery problem}, see e.g. \cite{MR}. 
This allows us to formulate optimality benchmarks
that steer our development of recovery algorithms. 
}
   
\subsection{Mathematical formulation of the state estimation problem}

The  {\em sensing} or {\em recovery} problem studied in this paper
are formulated in a Hilbert space $V$
equiped with some norm $\|\cdot\|$ and inner product $\<\cdot,\cdot\>$: we want to recover
an approximation to an unknown function $u\in V$
from data given by $m$ linear measurements
\be
\ell_i(u), \quad i=1,\dots,m,
\ee
where the $\ell_i$ are $m$ linearly independent linear functionals over $V$. 
This problem appears in many different settings. The particular one 
that motivates our work is the case where $u=u(y)$ represents the {\em state} of a physical system described
as a solution
to a parametric PDE 
\be
\cP(u,y)=0
\ee
for some unknown finite or infinite dimensional
parameter vector $y=(y_j)_{j\geq 1}$ picked 
from some admissible set $Y$. The $\ell_i$ are a mathematical
model for sensors that capture some partial information on the 
unknown solution $u(y)\in V$.

Denoting by $\omega_i\in V$ the Riesz
representers of the $\ell_i$, such that $\ell_i(v)=\<\omega_i,v\>$ for all $v\in V$, and defining
\be
W:={\rm span}\{\omega_1,\dots,\omega_m\},
\ee
the measurements are equivalently represented by
\be
w=P_{W}u. 
\ee
where $P_W$ is the orthogonal projection from $V$ onto $W$. A {\it recovery algorithm} is a computable map
\be
A: W \to V
\ee
and the approximation to $u$ obtained by this algorithm is
\be
u^*=A(w)=A(P_Wu).
\ee
The construction of $A$ should be based on the available prior 
information that describes the properties of the unknown $u$, and the evaluation of its performance needs
to be defined in some precise sense. Two distinct approaches are usually followed: 
\begin{itemize}
\item
In the {\it deterministic setting}, the sole prior information 
is that $u$ belongs to the set 
\be
\cM:=\{u(y)\; : \; y\in Y\},
\label{Mparam}
\ee 
of all possible solutions. The set $\cM$ is sometimes called the {\em solution manifold}. 
The performance of an algorithm $A$ over the class $\cM$ is usually measured by
the ``worst case'' reconstruction error
\be
E_{\wc}(A,\cM)=\sup\{\|u-A(P_Wu)\|\; : \; u\in\cM\}.
\ee
The problem of finding an algorithm that minimizes $E_{\wc}(A)$ is 
called {\it optimal recovery}. It has been extensively studied for convex sets $\cM$ that are 
balls of smoothness classes \cite{BB,MR,NW}, which is not the case for \iref{Mparam}.
\item
In the {\it stochastic setting}, the prior information on $u$ is 
described by a probability distribution $p$ on $V$,
which is supported on $\cM$, typically induced by a probability distribution on $Y$ that is assumed
to be known.
It is then natural to measure the performance of an algorithm
in an averaged sense, for example through the mean-square error
\be
E_{\ms}(A,p)=\E(\|u-A(P_Wu)\|^2)=\int_V \|u-A(P_W u)\|^2 dp(u).
\label{meansquare}
\ee
This stochastic setting is the starting point {\wnew for} {\it Bayesian estimation} 
methods \cite{DS}. Let us observe that for any algorithm $A$ one has $E_{\ms}(A,p)\leq E_{\wc}(A,\cM)^2$.
\end{itemize}

\subsection{Optimal algorithms}

The present paper concentrates on the deterministic setting according to the
above distinction, although some remarks will be given on the analogies
with the stochastic setting. In this setting, the benchmark for the performance
of recovery algorithms is given by
$$
E^*_{\wc}(\cM)=\inf_{A}E_{\wc}(A,\cM),
$$
where the infimum is taken over all possible maps $A$.

There is a simple mathematical description of an
optimal map that meets this benchmark.
For any bounded set $S\subset V$ we define
its {\it Chebychev ball} as the smallest closed ball that contains $S$. The 
{\it Chebychev radius and center} denoted by ${\rm rad}(S)$ and ${\rm cen}(S)$
are the radius and center of this ball. Since the information that we have
on $u$ is that it belongs to the set
\be
\cM_w:=\cM\cap V_w, \quad V_w:=\{v\in V\; : \; P_W v=w\}=w+W^\bot,
\ee
where $W^\perp$ is the orthogonal complement of $W$ in $V$, it follows that an optimal reconstruction map $A^*_\wc$ for the worst case error is given by
\be
A^*_\wc(w)={\rm cen}(\cM_w),
\ee 
because the Chebychev center of $\cM_w$ minimizes the quantity $\sup\{\|u-v\| \, :\, u\in\cM_w\}$ among all $v\in V$.
The worst case error is therefore given by
\be
E_\wc^*(\cM)=E_\wc(A^*_\wc,\cM)= \sup \{\rad(\cM_w)\, : \, w\in P_W (\cM)\}.
\ee
Note that the map $A^*_\wc$ is also optimal among all algorithms for each $\cM_w$, $w\in P_W (\cM)$, since 
\be
E_\wc(A^*_\wc,\cM_w)=\min_{A}E_\wc(A,\cM_w)=\rad(\cM_w), \quad w\in P_W(\cM).
\ee
However, there may exist other maps
$A$ such that $E_\wc(A,\cM)=E_\wc^*(\cM)$, since we 
also supremize over $w\in P_W(\cM)$.

\subsection{Linear and affine algorithms based on reduced models}

In practice the above map $A^*_\wc$ cannot be easily constructed
due to the fact that the solution manifold $\cM$ is a high-dimensional and geometrically 
complex object. One is therefore interested in designing ``sub-optimal yet
good'' recovery algorithms and analyze their performance. 

One vehicle for constructing linear recovery mappings $A$ is to use {\it  reduced modeling}. 
Generally speaking, 
reduced models consist of linear spaces $(V_n)_{n\geq 0}$ with
increasing dimension $\dim(V_n)=n$ which uniformly approximate
the solution manifold in the sense that
\be
\dist(\cM,V_n):=\max_{u\in \cM}\min_{v\in V_n}\|u-v\|\leq \e_n,
\ee
where 
\be
\e_0\geq \e_1\geq \cdots \geq \e_n\geq \cdots \geq 0,
\ee
are known tolerances. Instances of reduced models for parametrized families of PDEs
with provable accuracy are provided by polynomial approximations in the $y$ variable
\cite{CD,CDS1} or reduced bases \cite{BMPPT,S,RHP}. The construction of a reduced model
is typically done offline, using a large training set of instances of $u\in\cM$
called {\it snapshots}. The offline stage potentially has a high
computational cost. Once this is done, the online cost of recovering
$u^*=A(w)$ from any data $w$ using this reduced model should 
in contrast be moderate. 

In \cite{MPPY}, a simple reduced-model based recovery algorithm was proposed,
in terms of the map
\be
A_n(w):={\rm argmin}\{ \dist(v,V_n) \; :\; v\in V_w\},
\label{onespacelin}
\ee
which is well defined provided that $V_n\cap W^\perp =\{0\}$.  It turns out that $A_n$ is a linear mapping 
and so these algorithms are linear. {\wnew This approach is called the Parametrized-Background Data-Weak (PBDW)
method, however, we follow the terminology introduced in \cite{BCDDPW2},
refering to an algorithm of the form $A_n$ as {\it one-space-algorithm}.
In the latter,} it was shown that $A_n$ has a simple interpretation in terms of the cylinder 
\be
\cK_n:=\{ v\in V \; : \: {\rm dist}(v,V_n)\leq \e_n\},
\ee
that contains the solution manifold $\cM$. Namely, the algorithm $A_n$ is also given by
\be
A_n(w)=\cen(\cK_{n,w}), \quad \cK_{n,w}:= \cK_n\cap V_w,
\ee
{\mnew
and the map is shown to be the optimal when $\cM$ is replaced by the simpler
containement set $\cK_n$, that is
$$
A_n = \argmin_{A:W\to V} E_{\wc}(A, \cK_n). 
$$
}
The substantial advantage of this approach is that,
in contrast to $A^*_\wc$, the map $A_n$ can be easily computed by solving simple least-squares 
minimization problems which amount to finite linear systems. In turn $A_n$ is a linear 
map from $W$ to $V$. This map depends on $V_n$ and $W$, but not on $\e_n$ in view of \iref{onespacelin}. 
We refer to $A_n$ as the {\wnew{\it one-space-algorithm}}
based on the space $V_n$.

This algorithm satisfies the performance bound
\be
\|u-A_n(P_Wu)\|\leq \mu_n {\rm dist}(u,V_n\oplus (V_n^\perp \cap W))
\leq \mu_n {\rm dist}(u,V_n) \leq \mu_n\e_n,
\ee
where the last inequality holds when $u\in \cM$. Here 
\be
\mu_n=\mu(V_n,W):=\max_{v\in V_n} \frac{\|v\|}{\|P_W v\|},
\ee
is the inverse of the inf-sup constant $\beta_n:=\min_{v\in V_n} \max_{w\in W} \frac{\<v,w\>}{\|v\| \,\|w\|}$ 
which describes the angle between $V_n$ and $W$. In particular $\mu_n=\infty$ in the event where $V_n\cap W^\perp$
is non-trivial.

An important observation is that the {\wnew one-space} algorithm 
\iref{onespacelin} has a simple extension to the setting where  $V_n$ is an affine space  rather than a linear space, namely,  when 
\be
V_n=\o u + \wt V_n,
\label{affine}
\ee
with $\wt V_n$ a linear space of dimension $n$ and $\o u$ a given offset that is known to us.  

{\mnew
At a first sight, affine spaces do not bring any significant improvement in terms of approximating the solution manifold,
due to the following elementary observation: if $\cM$ is approximated with accuracy $\e$ by an $n$-dimensional affine space $V_n$ given by \iref{affine},
it is also approximated with accuracy $\wt \e\leq \e$ by the $n+1$-dimensional linear space
\be
\wt V_{n+1}:=V_n\oplus \R \o u.
\ee
However, the choice of an affine reduced model may significantly improve the performance of the 
{{\wnew one-space} algorithm} in the case where the parametric solution $u(y)$ is a ``small perturbation''
of a nominal solution $\o u=u(\o y)$ for some $\o y\in Y$, in the sense that 
\be
{\rm diam}(\cM)\ll \|u\|.
\ee
Indeed, suppose in addition that $\o u$ is badly aligned
with respect to the measurement space $W$ in the sense that 
\be
\|P_W \o u \| \ll \|u\|.
\ee
In such a case, any linear
space $V_n$ that is well tailored to approximating the solution manifold 
(for example a reduced basis space) will contain a direction
close to that of $\o u$ and thus, we will have that $\mu_n \gg 1$, rendering the
reconstruction by the linear one-space method much less accurate than
the approximation error by $V_n$. The use of the affine mapping
\iref{affine} has the advantage of elimitating the bad direction $\o u$ since
$\mu_n$ will now be computed with respect to the linear part $\wt V_n$.

A further perspective, currently under investigation, is to agglomerate local affine models in order
to generate 
{{{\em nonlinear}}} reduced model. This can be executed, for example, by
decomposing the parameter domain $Y$ into $K$ subdomains $Y_k$
and using different affine reduced models for approximating the resulting 
subsets $\cM_k=u(Y_k)$. 
}

\subsection{Objective and outline}

The standard constructions of reduced models are targeted at making the spaces
$V_n$ as efficient as possible for approximating $\cM$, that is, making $\e_n$ as small as possible
for each given $n$.
For example, for the reduced basis spaces, it is known \cite{BCDDPW1,DPW2}
that a certain greedy selection of snapshots
generates spaces $V_n$ such that $\dist(\cM,V_n)$ decays at the same rate (polynomial or exponential)
as the Kolmogorov $n$-width 
\be
\delta_n(\cM):=\inf\{\dist(\cM,E)\; : \; \dim(E)=n\}.
\label{kol}
\ee
However these constructions do not ensure the control of $\mu_n$ and therefore these reduced spaces may
be much less efficient when using the one-space algorithm for the recovery problem. 

In view of the above observations, the objective of this paper is to discuss
the construction of reduced models (both linear and affine) that are better targeted towards the 
recovery task. In other words, we want to build the spaces $V_n$ to make
the recovery algorithm $A_n$ as efficient as possible, given the measurement space $W$.
Note that a different problem is, given $\cM$, to optimize
the choice of the measurement functionals $\ell_i$ picked from
some admissible dictionary, which amounts to optimizing the space $W$,
as discussed for example in \cite{BCMN}. Here, we consider our 
measurement system to be imposed on us, and therefore $W$ to be fixed
once and for all. 

The rest of our paper is organized as follows.
In \S 2, we detail the affine map $A_n$ associated to $V_n$,
that can be computed in a similar way as in the linear case. Conversely, we show
that any affine recovery map may be interpreted as a one-space algorithm
for a certain affine reduced model $V_n$. 
For a general set $\cM$, the existence and construction of an optimal affine recovery map  $A^*_\wca$ for
the worst case error is therefore equivalent to the existence and construction of an optimal reduced space for the recovery problem. 
We then draw a short comparison with the stochastic setting in which the optimal
affine map $A^*_\msa$ for the mean-square error \iref{meansquare} is
derived explicitely from the second order statistics of $u$.

In \S 3, we compute an approximation of $A^*_\wca$ by convex optimization, 
based on a training set of snapshots. 
Two algorithms are considered: subgradient descent and primal-dual proximal splitting. 
Our numerical results illustrate the superiority of the
latter for this problem.  The optimal affine map $A^*_\wca$
significantly outperforms the one-space algorithm $A_{n^*}$
when standard reduced basis spaces $V_n$ are used and an optimal value $n^*$
is selected using the training set. 
It also outperforms the affine map $A_\msa^*$ computed 
from second order statistics of the training set. All three maps
significanly outperform the minimal $V$-norm recovery given by $A(w)=w=P_W u$.

\section{Affine {\wnew one-space} recovery}

 In this section, we show that any linear {\wnew recovery} algorithm is  given by a one-space algorithm and that a similar result holds for any affine algorithm.  We then go on to describe the optimal {\wnew one-space} algorithms by exploiting this
fact.

\subsection{The one-space algorithm}

We begin by discussing in more detail  the one-space algorithm for a   linear space $V_n$ of dimension $n\le m$.
As shown in \cite{BCDDPW2}, the map $A_n$  associated to $V_n$ has a simple expression after a proper choice of
favorable bases has been made for $W$ and $V_n$ through an SVD applied
to the cross-gramian of an initial pair of orthonormal bases. The resulting {\it favorable bases}
$\{\psi_1,\dots,\psi_m\}$  for $W$ and $\{\vp_1,\dots,\vp_n\}$ for $V_n$ satisfy the equations
\be
\<\psi_i,\vp_j\>=s_i\delta_{i,j},
\ee
where
\be
1\geq s_1\geq s_2\geq \dots \geq s_n>0 ,
\ee
are the singular values of the cross-gramian.
Then, if  $w$ is in $W$, we  can write $ w=\sum_{j=1}^m w_j\psi_j$ in the favorable basis,  and find that 
\be
A_n(w)= \sum_{j=1}^n s_j^{-1} w_j\vp_j+\sum_{j=n+1}^m w_j\psi_j.
\label{onespacefav1}
\ee
Let us observe that the functions $\psi_j$ in the second sum span the space $ V_n^\perp\cap W$ {\wnew while the first sum is the solution of the
least squares problem $\min_{v\in V_n} \|w- P_W v\|$ corrected by the second sum so as to fit the data.}


{\mnew
Now consider any linear recovery algorithm $A:W\to V$.
Since we are given the measurement observation $w$,
any algorithm $A$ which is a candidate to optimality must satisfy $P_W(A(w)) = w$
(otherwise the reconstruction error would not be minimized). Thus $A$ should have the form
\be
\label{opform}
A(w)=w+B(w),
\ee
where $B: W\to W^\perp$  with $W^\perp$ the orthogonal complement of $W$ in $V$. Note that in Functional Analysis the mappings $A$ of the form \eref{opform} are called liftings.

Therefore, in going further in this paper, we always require that $A$ has the form \eref{opform} and concentrate on the construction of good linear maps $B$.}
Our next observation is that any algorithm  $A$ of this form can always be interpreted
as a {\wnew one-space} algorithm  $A_n$ for a certain  space $V_n$ with $n\le m$. 

\begin{prop}
\label{proplin}
Let $A$ be any linear  map of the form  \eref{opform}.   Then, there exists 
a space $V_n$ of dimension $n\le m$ such that $A$ coincides with the one-space
algorithm \eref{onespacefav1} for $V_n$.
\end{prop}

\noindent
{\bf Proof:}
 By considering the SVD of the linear transform $B$, there exists an orthonormal basis
$\{\psi_1,\dots,\psi_m\}$ of $W$ and an orthonormal system $\{\omega_1,\dots,\omega_m\}$ in $W^\perp$
such that, with $w=\sum_{j=1}^m w_j \psi_j$,
\be
Bw=\sum_{j=1}^m\alpha_j w_j \omega_j, \quad  w\in W,
\ee
for some numbers $\alpha_1\geq \alpha_2 \geq \cdots \geq \alpha_m \geq 0$. Defining the functions
\be
\vp_j=s_j(\psi_j+\alpha_j\omega_j), \quad s_j=(1+\alpha_j^2)^{-1/2},
\ee
and defining $V_n$ as the span of those  $\vp_j$ for which  $\alpha_j\neq 0$, we recover the exact form \iref{onespacefav1}
of the one-space algorithm expressed in favorable bases.
\hfill $\Box$
\newline

These results can be readily extended to the case where $V_n$ is an affine space given by \iref{affine}
for some given $n$-dimensional linear space $\wt V_n$
and offset $\o u$. In what follows, we systematically use the notation
\be
\wt u=u-\o u,
\ee
for the recentered state, and likewise $\wt w=w-\o w$ with $\o w=P_W\o u$ for the recentered observation.
The one-space algorithm associated to $V_n$ has the form
\be
A_n(w)=\o u + \wt A_n(\wt w),
\label{affan}
\ee
where $\wt A_n$ is the one-space linear algorithm associated to $\wt V_n$. 

Performances bounds similar to those of the linear are derived 
in the same way as in \cite{BCDDPW2}: the reconstruction satisfies
\be
\|u-A_n(P_Wu)\|\leq \mu_n{\rm dist}(u,\o u + \wt V_n \oplus(\wt V_n^\perp\cap W))\leq \mu_n{\rm dist}(u,V_n),
\label{eq:err-one}
\ee
where 
\be
\mu_n=\mu(\wt V_n,W)=\max_{v\in \wt V_n} \frac{\|v\|}{\|P_W v\|}=s_n^{-1}<\infty.
\ee
The map $A_n$ is optimal for the cylinders of the form
\be
\cK_n=\{u\in V \; : \; {\rm dist}(u,V_n)\leq \e_n\},
\label{Kn}
\ee
since it coincides with the Chebychev
center of $\cK_{n,w}=\cK_n\cap V_w$. In particular, one has
\be
E^*_\wc(\cK_n)=E_\wc(A_n,\cK_n)=\mu_n\e_n.
\ee
For a solution manifold $\cM$ contained in $\cK_n$, one has
\be
E^*_\wc(\cM) \leq E_\wc(A_n,\cM)\leq \mu_n{\rm dist}(\cM, \wt V_n \oplus(\wt V_n^\perp\cap W))\leq \mu_n{\rm dist}(\cM, V_n)
\leq \mu_n\e_n,
\ee
and these inequalities are generally strict.

In view of \iref{affan} the map $A_n$ is affine. A general affine recovery map
takes form
\be
A(w)=w+Bw+c,
\label{genaffine}
\ee
where $B: W\to W^\perp$ is linear and $c=A(0)\in W^\perp$. The following result is a direct consequence
of Proposition \ref{proplin}.

\begin{cor}
\label{anycor}
Let $A$ be an affine map of the form \iref{genaffine}. Then, there exists 
an affine space $V_n=\o u+\wt V_n$ such that $A$ coincides with the one-space
algorithm \iref{affan}.
\end{cor}

\subsection{The best affine map}

In view of this result, the search for an affine reduced model $V_n$ that
is best tailored to the recovery problem is equivalent to the search of
an optimal affine map. Our next result is that such a map always 
exist when $\cM$ is a bounded set.

\begin{theorem}
Let $\cM$ be a bounded set. Then there exists a map $A^*_\wca$ that minimizes
$E_\wc(A,\cM)$ among all affine maps $A$.
\end{theorem}

\noindent
{\bf Proof:}  We consider any affine map $A$ of the form \iref{genaffine}, so that the error
is given by
\be
E_\wc(A,\cM)=\sup \{u\in \cM \; : \; \|P_{W^\perp}u-c-BP_W u\|\}=:F(c,B).
\ee
We begin by remarking that for each $(c,B)\in W^\perp \times \cL(W,W^\perp)$, the map
$u\mapsto \|P_{W^\perp}u-c-BP_W u\|$ is uniformly bounded 
on the bounded set $\cM$. Its supremum $F(c,B)$ is thus a finite 
positive number, which we may write as
\be
F(c,B)=\sup_{u\in\cM} F_u(c,B),
\label{F}
\ee
where $F_u(c,B)=\|P_{W^\perp}u-c-BP_W u\|$. 
Each $F_u$ is convex and satisfies the Lipschitz bound
\be
|F_u(c,B)-F_u(c',B')| \leq \|c-c'\| + M \|B-B'\|_{S},
\ee
with 
\be
\|B\|_S=\max\{\|Bv\| \; : \; v\in W, \; \|v\|=1\},
\ee 
the spectral norm and $M:=\sup\{\|P_W u\| \; : \; u\in\cM\}<\infty$.
This implies that the function $F$ is convex and satisfies the same Lipschitz bound.

We note that the linear maps of $\cL(W,W^\perp)$ are of rank at most $m$ and therefore, given any orthonormal basis $(e_1,\dots,e_m)$
of $W$, we can equip $\cL(W,W^\perp)$ with the Hilbert-Schmidt norm
\be
\|B\|_{HS}:=\(\sum_{i=1}^m \|Be_i\|^2\)^{1/2},
\ee
which is equivalent to the spectral norm since
\be
\|B\|_{S}\leq \|B\|_{HS}\leq \sqrt m \|B\|_{S}, \quad   B\in \cL(W,W^\perp).
\ee
In particular $F$ is continuous with respect to the Hilbertian norm 
\be
\|(c,B)\|_H:=\(\|c\|^2+\sum_{i=1}^m \|Be_i\|^2\)^{1/2}.
\ee
The function $F$ may not be infinite at infinity: this happens
if there exists a non-trivial pair $(c,B)$ such that
$$
c+BP_Wu=0, \quad u\in \cM.
$$
In order to fix this problem, we define
the subspace 
\be
S_0:=\Big\{(c,B)\in W^\perp\times \cL(W,W^\perp) \; : \; c+BP_Wu=0, \; u\in\cM\Big\}.
\ee
and we denote by $S_1$ its orthogonal complement in $W^\perp\times \cL(W,W^\perp)$
for the inner product associated to the above Hilbertian norm $\|\cdot\|_H$.
The function $F$ is constant in the direction
of $S_0$ and therefore we are left to prove the existence of the minimum
of $F$ on $S_1$. For any $(c,B)\in S_{1}$, there exists 
$u\in \cM$ such that 
$c+BP_Wu\neq 0$. This implies that
\be
\lim_{|t|\to +\infty} \|P_{W^\perp}u-tc-t BP_W u\| =+\infty,
\ee
and therefore that $\lim_{|t|\to +\infty}F_u(t(c,B))=+\infty$. This shows that $F$ is infinite
at infinity when restricted to $S_1$. Any convex and continuous function
in a Hilbert space is weakly lower semi-continuous, and admits a minimum when
it is infinite at infinity. We thus conclude in the existence
of a minimizer $(c^*,B^*)$ of $F$ and therefore 
\be
A^*_\wca (w)=w+c^*+B^* w,
\ee
is an optimal affine recovery map.
\hfill $\Box$

\subsection{\mnew The best affine map in the stochastic setting}

In the stochastic setting, assuming that $u$ has finite second order
moments, the optimal map that minimizes
the mean square error \iref{meansquare} is given by
the conditional expectation
\be
A^*_\ms(w)=\E( u \; | \; P_Wu=w),
\ee
that is, the expectation of posterior distribution $p_w$ of $u$
conditioned to the observation of $w$.
Various sampling strategies have been developed in order to approximate
the posterior and its expectation, see \cite{DS} for a survey. 
These approaches come at a significant computational cost {\wnew since}
they require a specific sampling for each instance 
$w$ of observed data. In the parametric PDE setting, each sample
requires one solve of the forward problem.

On the other hand, it is well known that an optimal
affine map $A^*_\msa$ for the mean square error can be explicitely derived from the first and second order statistics
of $u$. We briefly recall this derivation by using an arbitrary orthonormal basis $(e_1,\dots,e_m)$
of $W$ that we complement into an orthonormal basis $(e_j)_{j\geq 1}$ of $V$. We write
\be
u=\sum_{j\geq 1} w_j e_j \quad {\rm and}\quad
\o u=\E(u)=\sum_{j\geq 1} \o w_j e_j, \quad \o w_j:=\E(w_j),
\ee
as well as 
\be
\wt u=u-\o u=\sum_{j\geq 1} \wt w_j e_j, \quad \wt w_j:=w_j-\o w_j.
\ee
An affine recovery map of the form \iref{genaffine}
leaves the coordinates $w_1,\dots,w_m$ unchanged and recovers for each $i\geq 1$
\be
w_{m+i}^*=c_i+\sum_{j=1}^m b_{i,j} w_{j},
\ee
which can be rewritten as
\be
w_{m+i}^*=\o w_{m+i}+d_i+\sum_{j=1}^m b_{i,j} \wt w_{j}.
\ee
Since $E_\ms(A)=\sum_{i\geq 1}\E(|w_{m+i}^*-w_{m+i}|^2)$, the numbers $d_i$ and $b_{i,j}$
are found by separately minimizing each term. By Pythagoras theorem one has
\be
\E(|w_{m+i}^*-w_{m+i}|^2)=|d_i|^2 + \E\(|\sum_{j=1}^m b_{i,j} \wt w_{j}- \wt w_{m+i}|^2\),
\ee
which shows that we should take $d_i=0$. Minimizing the second term
leads to the orthogonal projection equations
\be
\sum_{j=1}^m b_{i,j} t_{j,l}=t_{m+j,l}, \quad l=1,\dots,m.
\ee
which involve the entries of the covariance matrix
\be
\bS:=(t_{i,j}), \quad t_{i,j}:=\E(  \wt w_{i} \wt w_{j}).
\ee
Therefore, with the block decomposition
\be
\bS=\left(\begin{array}{cc} \bS_{1,1}& \bS_{1,2} \\ \bS_{2,1}& \bS_{2,2}\end{array}\right),
\ee
corresponding to the splitting of rows and columns from $\{1,\dots,m\}$ and $\{m+1,m+2,\dots\}$, one obtains
that the matrix $\bB=(b_{i,j})$ that defines the optimal affine map
satisfies $\bS_{1,1} \bB^{\rm T}=\bS_{1,2}$ and therefore, 
\be
\bB=\bS_{2,1}\bS_{1,1}^{-1}
\label{linearest}
\ee
where we have used the symmetry of $\bS$. In other words,
\be
A^*_\msa(w)=w+P_{W^\perp} \o u+B \tilde w, 
\ee
where the linear transform $B\in \cL(W,W^\perp)$ is represented by the matrix $\bB$ in the basis $(e_j)_{j\geq 1}$.

The optimal affine recovery map $A^*_\msa$ agrees with the optimal map $A^*_\ms$
in the particular case where $u$ has Gaussian distribution, therefore entirely characterized by
its average $\o u$ and covariance matrix $\bS$. To see this, assume for simplicity that $V$ is
finite dimensional. The distribution of $\bu=(w_j)_{j\geq 1}$ has density proportional to $\exp(-\frac 1 2\<\bT\wt \bu,\wt\bu\>)$
where $\bT=\bS^{-1}$. We expand the quadratic form into
\be
\frac 1 2 \<\bT \wt\bu,\wt\bu\>=\frac 1 2 \<\bT_{1,1}\wt\bw,\wt\bw\>+\<\bT_{2,1}\wt\bw,\wt\bw_\perp\>+
\frac 1 2\<\bT_{2,2}\wt\bw_\perp,\wt\bw_\perp\>,
\ee
where $\wt \bw_\perp=(\wt w_{m+j})_{j\geq 1}$ and
$\wt \bw=(\wt w_{j})_{j=1,\dots,m}$, and where
\be
\bT=\left(\begin{array}{cc} \bT_{1,1}& \bT_{1,2} \\ \bT_{2,1}& \bT_{2,2}\end{array}\right),
\ee
is a block decomposition similar to that of $\bS$. The distribution of 
the vector $\wt \bw_\perp$ conditional to the observation of $\wt \bw$
is also gaussian and its expectation
coincides with the minimum of the quadratic form
\be
Q_{\bw}(\wt \bw_\perp)=\frac 1 2\<\bT_{2,2}\wt\bw_\perp,\wt\bw_\perp\>+\<\bT_{2,1}\wt\bw,\wt\bw_\perp\>.
\ee
Therefore 
\be
\E(\wt\bw_\perp\,|\,\wt \bw)=-\bT_{2,2}^{-1}\bT_{2,1} \wt \bw=\bS_{2,1}\bS_{1,1}^{-1}\wt \bw=\bB \wt \bw,
\label{condest}
\ee
which shows that 
\be
A^*_\ms(w)=\E(u\,|\, P_W u=w)=A_\msa^*(w).
\ee

{\mnew One main interest of the above discussed stochastic setting is that the
best affine map is now explicitely given 
by the second order statistics, in view of \iref{linearest}.
This contrasts with the deterministic setting
in which the optimal affine map is obtained by
minimization of the convex functional $F$  from \iref{F}
and does not generally have a
simple explicit expression. Algorithms for solving this 
minimization problem are the object of
the next section.

Only for particular cases where $\cM$ has a simple geometry, the best affine map  $A^*_\wca$ in the
deterministic setting has a simple expression.
One typical example is when $\cM$ is an ellipsoid described
by an equation of the form
\be
\<\bT \wt\bu,\wt\bu\>\leq 1,
\ee
for a symmetric positive matrix $\bT$. Then, the set $\cM_w=\cM \cap V_w$ is also an ellipsoid associated
with the above quadratic form $Q_\bw$. The coordinates of its center are therefore given by
the equation $\wt\bw_\perp=-\bT_{2,2}^{-1}\bT_{2,1} \wt \bw$, which is 
the same as that defining the conditional expectation in \iref{condest}
This shows that, in the particular case of an ellipsoid, 
(i) the optimal map $A^*_\wc$ agrees with the optimal affine recovery map $A^*_\wca$ 
for the worst case error, and (ii) it has an explicit expression
which agrees with the optimal map $A^*_\msa$ for the mean square error
when the prior is a Gaussian with $\bT$ as inverse covariance matrix.
}

\section{Algorithms for optimal affine recovery}

\subsection{Discretization and truncation}
\label{sec:trunc}
We have seen that the optimal affine recovery map is obtained by minimizing the
convex function
\be
F(c,B)=\sup_{u\in\cM} \|P_{W^\perp}u-c-BP_W u\|,
\ee
over $W^\perp \times \cL(W,W^\perp)$. This optimization problem
cannot be solved exactly for two reasons:
\begin{enumerate}
\item
The sets $W^\perp$ as well as $\cL(W,W^\perp)$ are infinite dimensional
when $V$ is infinite dimensional.
\item
One single evaluation of $F(c,B)$ requires in principle to explore the entire manifold $\cM$.
\end{enumerate}

The first difficulty is solved by replacing $V$ by a subspace 
$Z_N$ of finite dimension $\dim(Z_N)=N$ that approximates the solution manifold $\cM$ with an accuracy
of smaller order than that expected for the recovery error. One possibility is
to use a finite element space $Z_N=V_h$ of sufficiently small mesh size $h$. However its resulting dimension $N=N(h)$
needed to reach the accuracy could still be quite large. An alternative is to
use reduced model spaces $Z_N$ which are more
efficient for the approximation of $\cM$, as we discuss further.

We therefore minimize $F(c,B)$ over $\wt W^\perp \times \cL(W,\wt W^\perp)$,
where $\wt W^\perp$ is the orthogonal complement of $W$ in the space $W+Z_N$,
and obtain an affine map $\wt A_\wca$ defined by
\be
\label{redopt}
\tilde  A_\wca(w)=w+\o c+\o Bw, \quad (\o c,\o B):=\argmin\{F(c,B): \ c\in \wt W^\perp, B\in\cL(W,\wt W^\perp)\}.
\ee

 In order to compare its performance with that of $A^*_\wca$,
we first observe that
\be
\label{eq:epsN}
\|P_{W^\perp}u-P_{\wt W^\perp}u\| \leq \e_N:=\sup_{u\in\cM}{\rm dist}(u,Z_N).
\ee
For any $(c,B)\in W^\perp \times \cL(W,W^\perp)$, we define
$(\wt c,\wt B)\in \wt W^\perp \times \cL(W,\wt W^\perp)$ by $\wt c=P_{\wt W^\perp}c$
and $\wt B=P_{\wt W^\perp}\circ B$. Then, for any $u\in \cM$,
$$
\begin{array}{ll}
\|P_{W^\perp}u-\wt c-\wt B u\| &\leq \|P_{\wt W^\perp}(P_{W^\perp}u-c-BP_W u)\|
+\|P_{W^\perp}u-P_{\wt W^\perp}u\| \\
&\leq  \|P_{W^\perp}u-c-BP_W u\|+\e_N.
\end{array}
$$
It follows that we have the framing
\be
E(A^*_\wca,\cM)\leq E(\wt A_\wca,\cM)\leq E(A_\wca^*,\cM)+\e_N,
\label{frame}
\ee
which shows that the loss in the recovery error is at most of the order $\e_N$.

To understand how large $N$ should be, let us observe that a recovery map
$A$ of the form \iref{genaffine} takes it value in the
linear space
\be
F_{m+1}=\R c + {\rm ran}(B),
\ee
which has dimension $m+1$. It follows that
the recovery error is always larger than 
the approximation error by such a space. Therefore
\be
E_{wc}(A^*_\wca,\cM) \geq \delta_{m+1}(\cM),
\ee
where $\delta_{m+1}(\cM)$ is the Kolmogorov $n$-width defined by \iref{kol} for $n=m+1$.
Therefore, if we could use the space $Z_n:=E_n$ that exactly achieve  
the infimum in \iref{kol}, we would be ensured that, with $N=m+1$, the 
additional error $\e_N=\delta_{m+1}(\cM)$ in \iref{frame} is of smaller order than 
$E_{wc}(A^*_\wca,\cM)$. As a result we would obtain the framing
\be
E(A^*_\wca,\cM)\leq E(\wt A_\wca,\cM)\leq 2E(A_\wca^*,\cM),
\label{frame2}
\ee
 In practice, since we do not have access to 
the $n$-width spaces, we use instead the reduced basis spaces $Z_n:=V_n$
which are expected to have comparable approximation performances in view of the results
from \cite{BCDDPW1,DPW2}. We take $N$ larger than $m$ but of comparable order.

The second difficulty is solved by replacing the set $\cM$ 
in the supremum that defines $F(c,B)$ by a discrete training set $\wt \cM$,
which corresponds to a discretization $\wt Y$ of the parameter domain $Y$, that is
\be
\wt \cM:=\{u(y)\, : \, y\in\wt Y\},
\ee
with finite cardinality.

We therefore minimize over $\wt W^\perp \times \cL(W,\wt W^\perp)$ the function
\be
\wt F(c,B)=\sup_{u\in\wt\cM} \|P_{W^\perp}u-c-BP_W u\|,
\ee
which is computable. The additional error resulting from this discretization can 
be controlled from the resolution of the discretization. Namely, let $\e>0$ be the smallest
value such that $\wt \cM$ is
an $\e$-approximation net of $\cM$, that is, $\cM$ is covered by the $V$-balls $B(u,\e)$
for $u\in \wt\cM$. Then, we find that
\be
\wt F(c,B)\leq F(c,B)\leq \wt F(c,B)+\e \|B\|_{\cL(W,\wt W^\perp)},
\ee
which shows that the additional recovery error will be of the order of $\e$
amplified by the norm of the linear part of the optimal recovery map.

One difficulty is that the cardinality of $\e$-approximation nets become potentially
untractable for small $\e$ as the parameter dimension becomes large, due to
the curse of dimensionality. This difficulty also occurs when
constructing reduced basis by a greedy selection process which also needs
to be performed in a sufficiently dense discretized sets. 
Recent results obtained in \cite{CDD} show 
that in certain relevant instances $\e$ approximation nets can
be replaced by random training sets of smaller cardinality. One interesting
direction for further research is to apply similar ideas in the context of the present paper.

\subsection{Optimization algorithms}
\label{sec:optim}
\begin{figure}[H]
\centering
\includegraphics[scale=0.5]{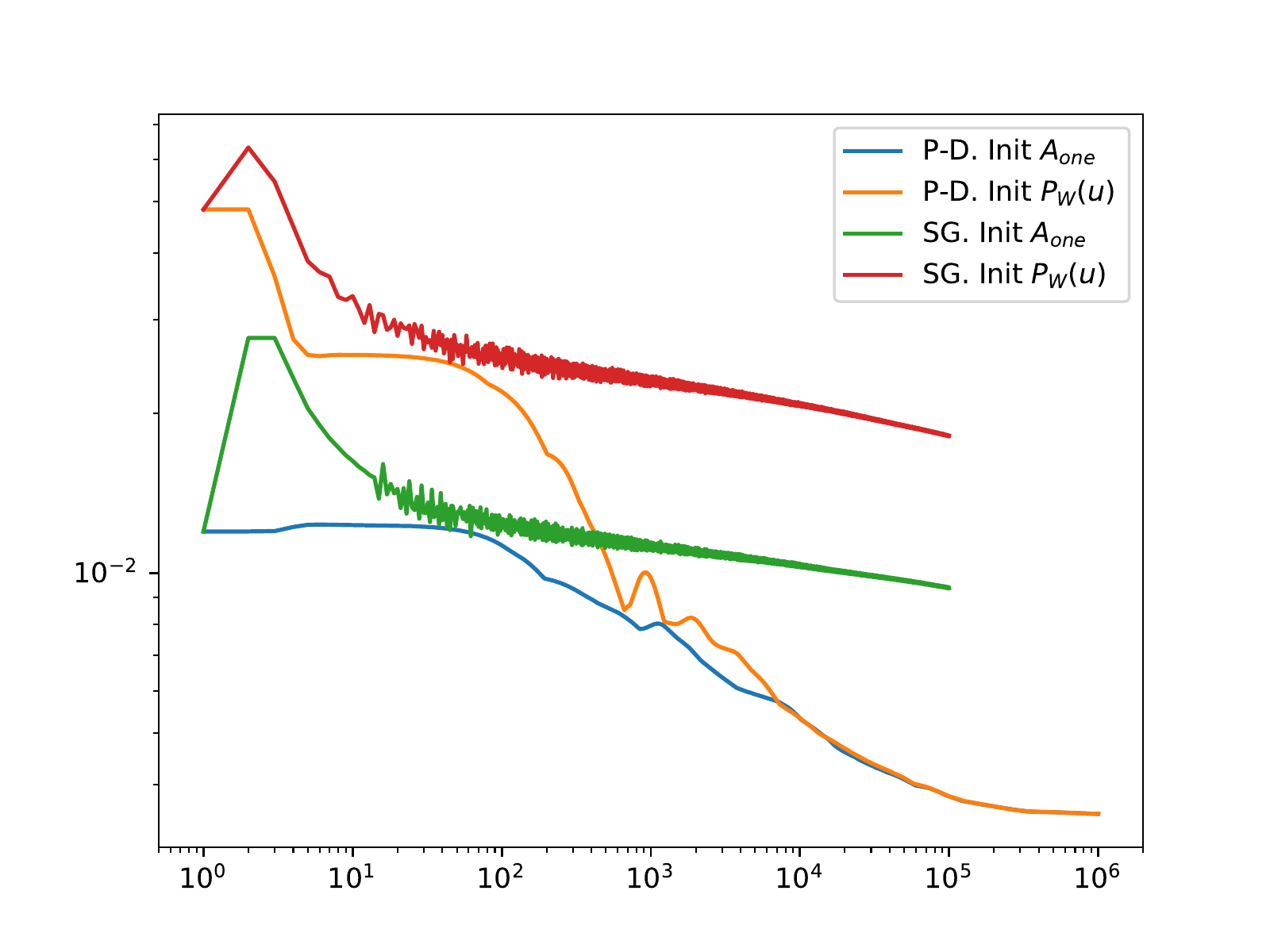}
\caption{Convergence of the objective function for two different optimization algorithms and starting guesses. P.D. = Primal-Dual splitting. S.G.=Subgradient. Here, $m=40$.}
\label{fig:conv-pd-starting-guess}
\end{figure}

As already brought up in the previous section, the practical computation of $\wt A_{\wc}$ consists in solving
\be\label{eq:opt-hilbert}
\min_{(c,B) \in \wt W^\perp \times \cL(W,\wt W^\perp)} \sup_{u\in\wt\cM} \|P_{W^\perp}u-c-BP_W u\|^2,
\ee
The numerical solution of this problem is challenging due to its lack of smoothness (the objective function is convex but non differentiable) and its high dimensionality (for a given target accuracy $\e_N$, the cardinality of $\wt \cM$ might be large). One could use classical subgradient methods, which are simple to implement.
However these schemes only guarantee a very slow $O(k^{-1/2})$ convergence rate of the objective function, where $k$ is the number of iterations. This approach did not give satisfactory results in our case: due to the slow convergence, the solution update of one iteration falls below machine precision before approaching the minimum close enough, see Figure \ref{fig:conv-pd-starting-guess}. This has motivated the use of a primal-dual splitting method which is known to ensure a $O(1/k)$ convergence rate on the partial duality gap. We next describe this method, but only briefly, as a detailed analysis would make us deviate too far from the main topic of this paper. A complete analysis with further examples of application will be presented in a forthcoming work \cite{FM2019}.

We assume without loss of generality that $\dim(W+V_N) = m+N$ and that $\dim \wt W^\perp = N$. Let $\{\psi_i\}_{i=1}^{m+N}$ be an orthonormal basis of $W+V_N$ such that $W={\rm span}\{\psi_1,\dots,\psi_m\}$. Since for any $u\in V$,
$$
P_{W+V_N} u = \sum_{i=1}^{m+N} u_i \psi_i,
$$
the components of $u$ in $W$ can be given in terms of the vector $\bw= ( u_{i} )_{i=1}^{m}$ and the ones in $\wt W^\perp$ with $\bu = ( u_{i+m} )_{i=1}^{N}$. 

We now consider the finite training set
\be
\wt\cM:=\{u^1,\dots,u^J\}, \quad J:=\#(\wt\cM)<\infty,
\ee
and denote by $\bw^j$ and $\bu^{j}$ the vectors
associated to the snapshot functions $u^j$ for $j=1,\dots,J$. 
One may express the problem \eqref{eq:opt-hilbert} as the search for
\begin{equation}
\label{eq:opt-Rn}
\min_{ \substack{(\bR, \bb) \in \\ \R^{N\times m}\times \R^{N}}} \max_{j=1,\dots,J} \Vert \bu^{j} - \bR \bw^{j} - \bb \Vert^2_{2}.
\end{equation}
Concatenating the matrix and vector variables $(\bR, \bb)$ into a single 
$\bx \in \R^{m(N+1)}$, we rewrite the above problem as
\begin{equation}
\label{eq:minP}
\min_{\bx \in \R^{m(N+1)}} \max_{j=1,\dots,J} f_j(\bQ_j\bx),
\end{equation}
where $\bQ_j\in \R^{N\times m(N+1)}$ is a sparse matrix built using the coefficients of $\bw^j$ and 
$f_j(\by) := \Vert \bu^{j} - \by \Vert^2_{2}$.

The key observation to build our algorithm is that problem~\eqref{eq:minP} can be equivalently written as a minimization problem on the epigraphs, i.e.,
\be
\begin{aligned}
	&\min_{(\bx,t) \in \R^{m (N+1)} \times \R^+} t \quad \text{subject to} \quad f_j(\bQ_j \bx) \leq t , \quad j=1,\dots,J \\
\iff&\min_{(\bx,t) \in \R^{m (N+1)} \times \R^+} t \quad \text{subject to} \quad (\bQ_j \bx,t) \in \epi_{f_j}, \quad j=1,\dots,J,
\end{aligned}
\ee
or, in a more compact (and implicit) form,
\begin{equation}\label{eq:minPepi}\tag{$\mathrm{P_{epi}}$}
\min_{(\bx,t) \in \R^{m (N+1)} \times \R^+} t + \sum_{j=1}^J \iota_{\epi_{f_j}}\pa{\bQ_j\bx,t} .
\end{equation}
where, for any non-empty set $S$ the indicator function $\iota_S$ has value $0$ on $S$ and $+\infty$ on $S^c$.
 
This problem takes the following canonical expression, which is amenable to a primal-dual proximal splitting algorithm
\begin{equation}\label{eq:primal-problem}
\min_{(\bx,t) \in \R^{m (N+1)} \times \R} G(\bx,t) + F \circ L(\bx,t).
\end{equation}
Here, $G$ is the projection map for the second variable
\be
G(\bx,t) = t,
\ee 
the linear operator $L$ is defined by
\be
L (\bx,t):=\pa{(\bQ_1 \bx,t),(\bQ_2 \bx,t),\cdots,(\bQ_J \bx,t)}
\ee
and acts from $\R^{m (N+1)} \times \R$ to $\times_{j=1}^J (\R^N \times \R)$
and the function $F$ acting from $\times_{j=1}^J (\R^N \times \R)$ to $\R$ is defined by
\be
F\( (\bv_1,t_1),\cdots,(\bv_J,t_J)\):= \sum_{j=1}^J \iota_{\epi_{f_j}}\pa{\bv_j,t_j}.
\ee
Note that $F$ is the indicator function of the cartesian product of epigraphs.

Before introducing the primal-dual algorithm, some remarks are in order:
\begin{enumerate}
\item We recall that if $\phi$ is a proper closed convex function on $\R^d$, its proximal mapping $\prox_\phi$ is defined by
\be
\prox_\phi(y)={\rm argmin}_{\R^d}\(\phi(x)+\frac 1 2\|x-y\|_2^2\).
\ee
\item The adjoint operator $L^*$ is given by
\be
L^*\( (\bv_1,t_1),\cdots,(\bv_J,t_J)\):= \pa{\sum_{j=1}^J \bQ_j^T \bv_j,\sum_{j=1}^J t_j} .
\ee
It can be easily shown that the operator norm of $L$ satisfies $\norm{L}^2 \leq J + \sum_{j=1}^J \norm{\bQ_j}^2$.
\item Both $G$ and $F$ are simple functions in the sense that their proximal mappings, $\prox_{G}$ and $\prox_{F}$, can be computed in closed form. See \cite{FM2019} for details.
\end{enumerate}


The iterations of our primal-dual splitting method read for $k\geq 0$,
\be
\label{eq:its-pd}
\begin{aligned}
(\bx,t)^{k+1} &= \prox_{\gamma_G G} \( (\bx, t)^k - \gamma_G L^* \( \( (\bv_1,\xi_1),\dots, (\bv_J,\xi_J) \)^k\) \) , \\
(\bar\bx,\bar t)^{k+1}  &= (\bx,t)^{k+1} + \theta \( (\bx,t)^{k+1} - (\bx,t)^k \) , \\
\( (\bv_1,\xi_1),\dots, (\bv_J,\xi_J)\)^{k+1} &= 
\prox_{\gamma_F \hat F}  
\(  \((\bv_1,\xi_1),\dots, (\bv_J,\xi_J)\)^k 
+ \gamma_F L (\bar\bx,\bar t)^{k+1} \), 
\end{aligned}
\ee
where $\hat F$ is the Fenchel-Legendre transform of $F$, $\gamma_G > 0$ and $\gamma_F > 0$ are such that $\gamma_G\gamma_F < 1/\norm{L}^2$, and $\theta \in [-1,+\infty[$ (it is generally set to $\theta=1$ as in \cite{chambolle2011first}).

{\mnew
Algorithm \ref{alg:pd} gives some guidelines and summarizes in an informal pseudo-code style the main steps of the primal-dual approach (the implementation of the routine ``BuildQ'' is left to the reader).}

\begin{algorithm}[H]
   \caption{Primal-dual algorithm: $\bR, \bb = \textsc{PD}(\widetilde \cM, \cM_{\greedy}, W, K_{\max})$}
  \label{alg:pd}
  \begin{algorithmic}[1]
   \STATE \textbf{Input:}
\begin{itemize}
\item training manifold $\widetilde \cM$ for primal dual iterations
\item training manifold $\cM_{\greedy}$ for greedy algorithm
\item basis $\{\omega_i\}_{i=1}^m$ of measurement space $W$
\item maximum number of iteration $K_{\max}$
\end{itemize}
   \STATE Generate basis $\{v_i\}_{i=1}^N$ of $V_N$ \COMMENT{e.g.~with a greedy algorithm over $\cM_{\greedy}$, see \eqref{eq:greedyAlg}}
   \STATE Build orthonormal basis $\{ \psi_i\}_{i=1}^{m+N} $ of $W+V_N$ with a Gram-Schmidt procedure over $\{w_1,\dots,w_m,v_1,\dots, v_N\}$. In this way, $\widetilde W^\perp = {\rm span} \{ \psi_{m+1},\dots, \psi_{m+N} \}$.
   \STATE Qlist, wlist, ulist = []
   \FORALL{$u \in \widetilde \cM$} \COMMENT{Build matrices $\bQ^j$ of \eqref{eq:minP}}
   	\STATE $\bw =  \{ \< u, \psi_i \>  \}_{i=1}^m$, $\bu =  \{ \< u, \psi_i \>  \}_{i=m+1}^{N+m}$, $\bQ = \mathrm{BuildQ}(\bw)$
	\STATE Qlist.append($\bQ$), wlist.append(\bw), ulist.append(\bu)
   \ENDFOR
   \State Estimate $\Vert L \Vert$ \COMMENT{e.g.~with power method}
   \State Set $\gamma_G$ and $\gamma_F$ such that $\gamma_G\gamma_F < 1 / \Vert L \Vert^2$
   \State $\bar \bx = \bx =$ zeros($m(N+1)$)  \COMMENT{starting guess for $A:W\to V$ set to $A(w)= w$.}
   \State $t=1$ \COMMENT{starting guess for $t>0$.}
   \State $\( (\bv_1,\xi_1),\dots, (\bv_J,\xi_J) \) = L(\bx, t)$ \COMMENT{starting guess dual variables.}
   \FOR{$k$ in $[0, K_{\max}]$} \COMMENT{primal-dual iterations}
   	\STATE $(\bx_{\text{old}}, t_{\text{old}}) = (\bx, t)$
   	\STATE $\( (\bv_1,\xi_1),\dots, (\bv_J,\xi_J)\)= \prox_{\gamma_F \hat F}  \(  \((\bv_1,\xi_1),\dots, (\bv_J,\xi_J)\) + \gamma_F L (\bar\bx,\bar t) \)$
	\STATE $(\bx,t) = \prox_{\gamma_G G} \( (\bx, t) - \gamma_G L^* \( \( (\bv_1,\xi_1),\dots, (\bv_J,\xi_J) \)\) \) $
	\STATE $(\bar\bx,\bar t) = (\bx,t) + \theta \( (\bx,t) - (\bx_{\text{old}}, t_{\text{old}}) \)$
   \ENDFOR
   \STATE Retrieve $\bR, \bc$ by appropriately reshaping \bx
   \STATE \textbf{Output:} $\bR, \bc$
  \end{algorithmic}
\end{algorithm}

To illustrate the relevance of this algorithm for our purposes, we compare its performance with a standard subgradient method. Figure \ref{fig:conv-pd-starting-guess} plots the convergence history of the objective function across the iterations of both optimization methods in the example described in the next section ($m=40$, $N=110$ and $J=10^3$). Two different reconstruction maps have been considered as starting guesses: 
the minimal $V$-norm recovery map given by $A(w)=w=P_W u$, and the 
one-space algorithm $A_{n^*}$ based on reduced basis spaces $V_n$ 
with an optimal choice $n^*$ for $n$. The convergence plot shows the superiority of the primal-dual method
which converges to the same minimal value of the objective function after $10^5$ iterations regardless of the intialization, while
the subgradient method fails to reach it since its increments fall below machine precision.

\begin{figure}[H]
\centering
\includegraphics[scale=0.5]{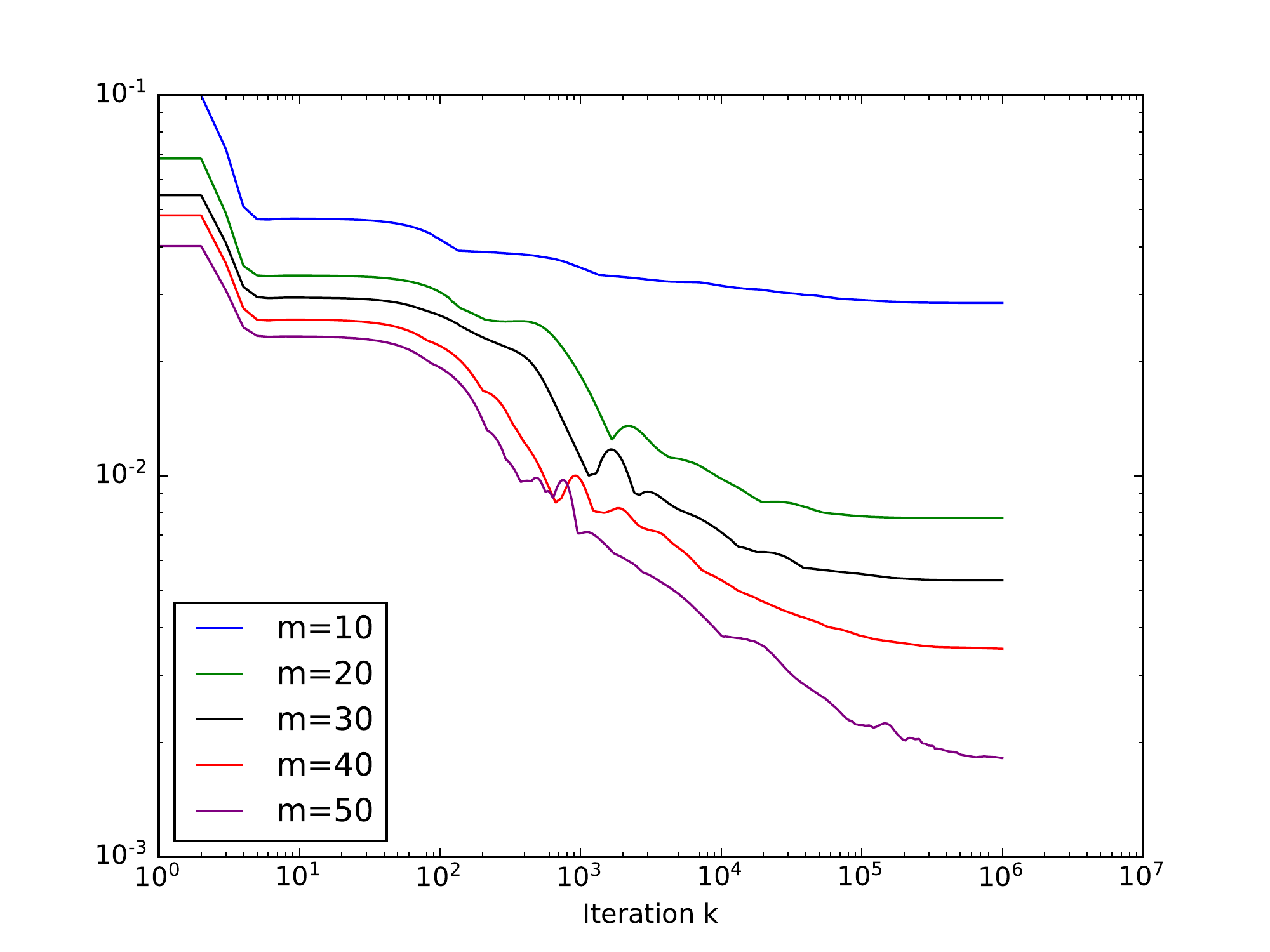}
\caption{Convergence of the objective function in the primal-dual iterations for $m=10,\,20,\,30,\,40,\,50$.}
\label{fig:conv-pd}
\end{figure}

For the same numerical example described next, we vary $m$ and consider $m=10,\,20,\,30,\,40,\,50$. Figure \ref{fig:conv-pd} gives the convergence of the reconstruction error over the training set  $\wt\cM$ across the primal-dual iterations (for simplicity, we took $P_{W_m}$ as the starting guess for $A^{(m)}_{\wca}$). To make sure that we reach convergence, we perform $10^6$ iterations for each case. As expected, we observe in this figure that the final value of the objective function decreases as we increase the value of $m$ (the reconstruction error decreases as we increase the number of measurements).

\subsection{Numerical tests}

We present some numerical experiments, aiming primarily at comparing in terms of the maximum reconstruction error the three above discussed recovery maps: the one-space affine map $A_n$, the best affine map $A^*_{\msa}$ for the mean-square error, and the best affine map $A^*_{\wca}$. for the 
worst case error. In addition, we also consider the mimimum $V$ norm reconstruction map $A(w)=w=P_Wu$. 
The results highlight the superiority of the best affine algorithm {\mnew with respect to the reconstruction error. This comes however} at the cost of a computationally intensive training phase as previously described.

We consider the elliptic problem
\be \label{diff_eq_y}
\begin{cases}
-{\rm div}\(a(y) \nabla u\)&=f, \quad  x\in D\\
\hfill u(x)&=0, \quad  x \in  \partial D
\end{cases}
\ee
on the unit square $D=]0,1[^2$, with a certain parameter dependence in the field $a$. More precisely, for a given $p\geq1$, we consider ``checkerboard'' random fields where $a(y)$ is piecewise constant on a $p\times p$ subdivision of the unit-square.
$$
D = \bigcup_{i,j=0}^{p-1} S_{i,j},
$$
with
$$
S_{i,j} := \Big[ \frac i p \, , \frac{i+1}{p} \Big[ \,\times\,
\Big[ \frac j p \, , \frac{j+1}{p} \Big[,\qquad i,j \in 0,\ldots,p-1.
$$
The random field is defined as
\be
a(y) = 1 + \frac{1}{2} \sum_{i,j = 0}^{p-1} \Chi_{S_{i,j}} y_{i,j},
\ee
where $\Chi_{S}$ denotes the characteristic function of a set $S$, and the $y_{i,j}$  are random coefficients that are independent, each with identical uniform distribution on $[-1,1]$. Thus, our vector of parameters is
$$
\by = (y_{i,j})_{i,j =0}^{p-1} \in \R^{p\times p}.
$$

In our numerical tests, we take $p=4$, that is $16$ parameters, and work in the ambient space  $V=H^1_0(D)$. All the sets of snapshots used for training and validating the reconstruction algorithms have been computed by first generating a certain number $J$ of random parameters $\by^{1},\ldots, \by^{J}$, with each $\by^{i} \in  [-1,1]^{p\times p}$, and then solving the variational form of \eqref{diff_eq_y} in $V=H^1_0(D)$ using $\mathbb{P}_1$ finite elements on a regular grid of mesh size $h=2^{-7}$. This gives the corresponding solutions $u_h^i = u_h(\by^{i})$ that are used in the computations. To ease the reading, in the following we drop the dependence on $h$ in the notation.

\begin{figure}[H]
\centering
\includegraphics[width=0.45\textwidth]{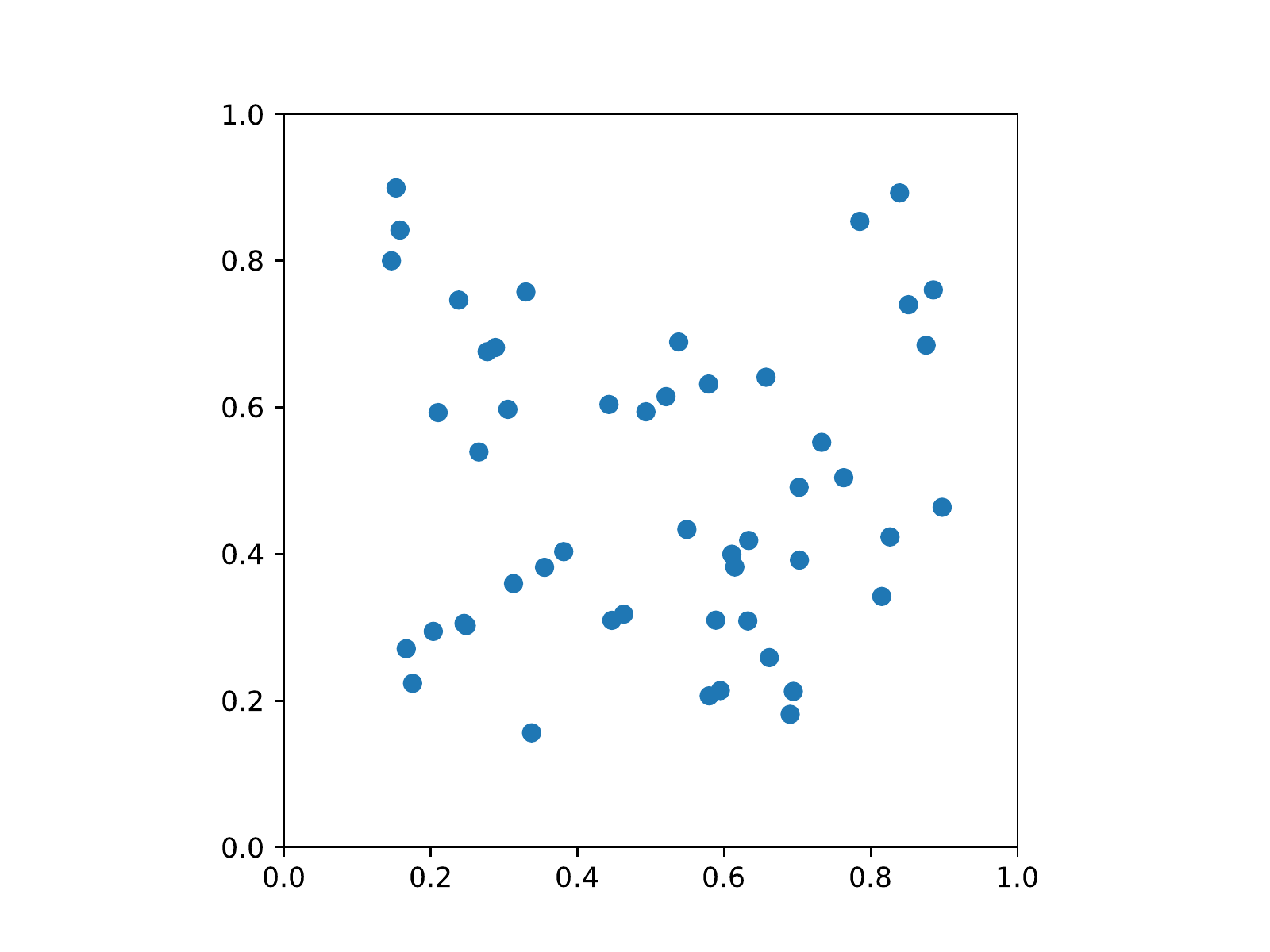}
\includegraphics[width=0.45\textwidth]{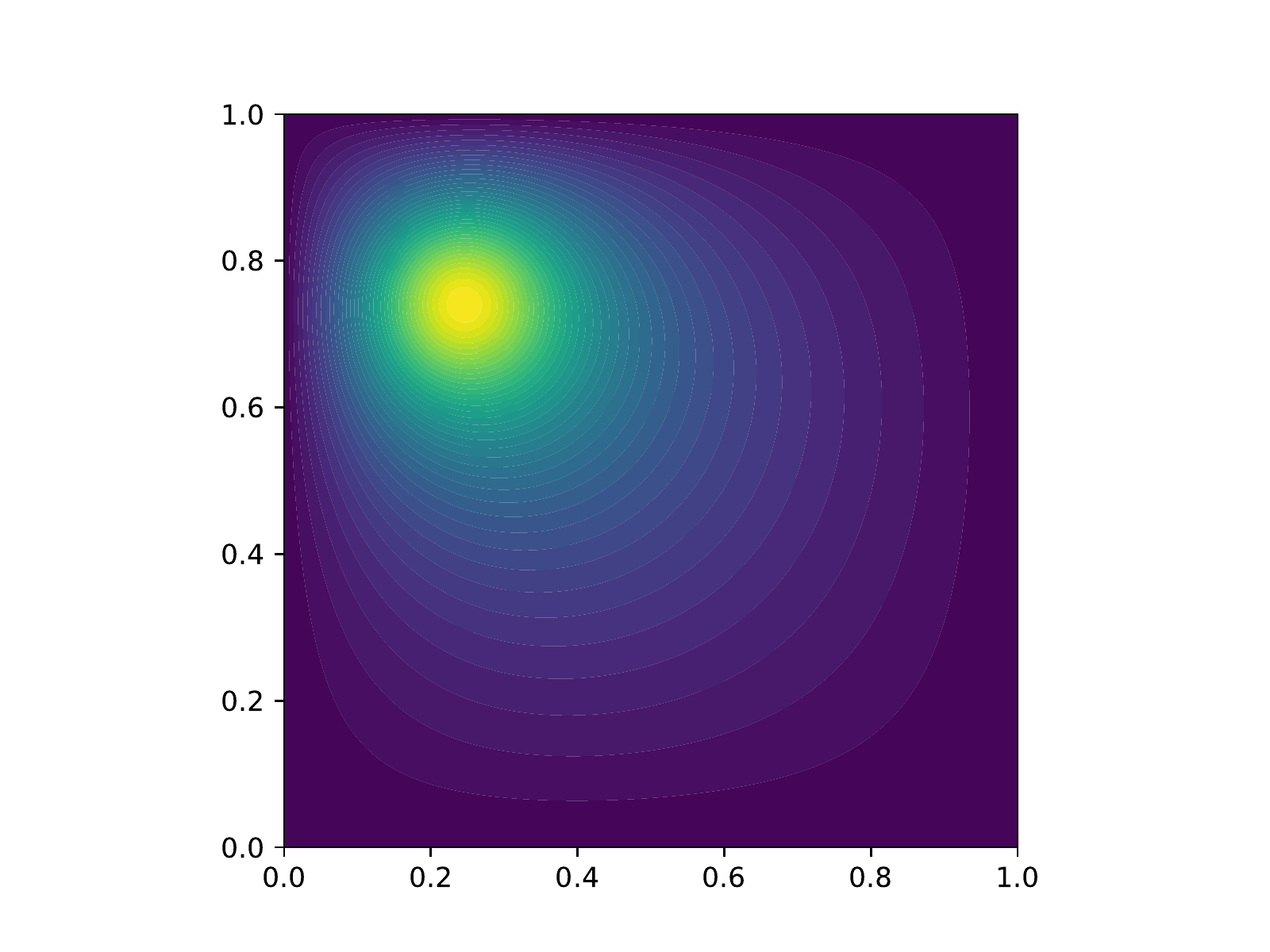}
\caption{Sensor locations and the function $\omega_{\bx_i,\tau_i}$ for $i=10$ ($\bx_i=(0.23, 0.75)$ and $\tau_i=0.06$).}
\label{fig:loc-sensors}
\end{figure}

The sensor measurements are modelled with linear functionals that are local averages of the form
\be
\ell_{\bx,\tau}(u)=\int_D u(\br) \vp_\tau(\br-\bx) \,\rm{d}\br,
\label{locav}
\ee
where 
\be
\vp_\tau(\br) \propto \exp( -| \br  |/2 \tau^2 )
\ee
is a radial function such that  $\int \vp_\tau=1$. The parameter $\tau>0$ represents the spread around the center $\bx$. For the observation space $W$ of our example, we randomly select $m=50$ centers $\bx_i \in [0.1, 0.9]^2$ and spreads $\tau_i\in [0.05, 0.1]$, and compute the Riesz representers $\omega_{\bx_i,\tau}$ of $\ell_{\bx_i,\tau}$ in $H^1_0(D)$. We then set
$$
W \coloneqq \{ \omega_{\bx_i,\tau} \}_{i=1}^M
$$
which is a space of dimension $m=50$. Figure \ref{fig:loc-sensors} shows the $m$ centers $\bx_i$. As an example, the figure also plots the function $\omega_{\bx_i,\tau}$ for $i=10$, which has center $\bx_i=(0.23, 0.75)$ and spread $\tau_i=0.06$.

\begin{figure}
\centering
\includegraphics[scale=0.5]{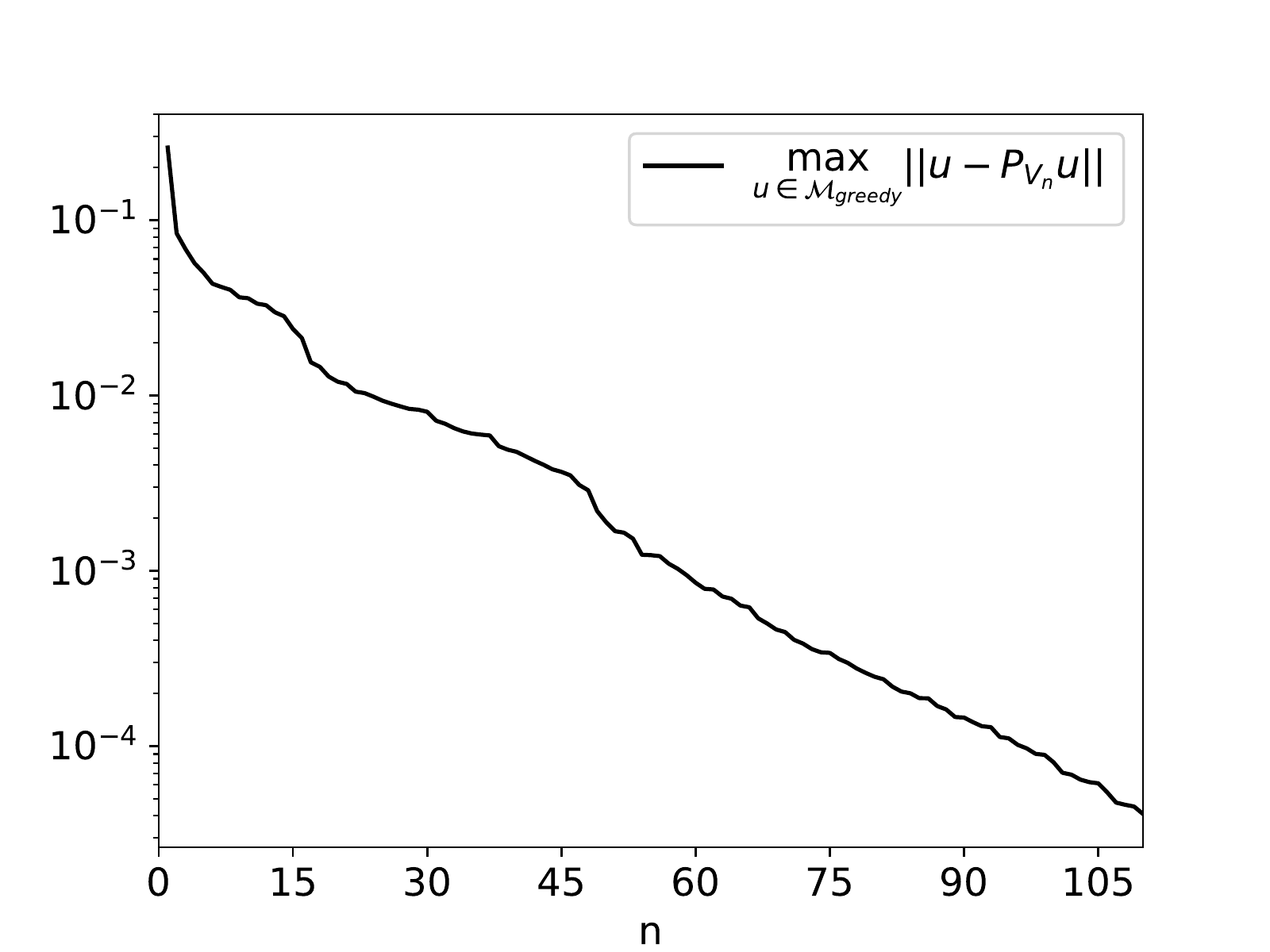}
\caption{Greedy algorithm: decay of the error $e_n^{\text{(greedy)}}=\max_{u\in\cM_{\greedy}} || u - P_{V_n}u ||$.}
\label{fig:greedy}
\end{figure}

As explained in section \ref{sec:trunc}, the first step to compute the best algorithm in practice consists in replacing $V=H^1_0(D)$ by a finite dimensional space that approximates the solution manifold $\cM$ at an accuracy smaller than the one expected for the recovery error. Here, we replace $V$ by $W+V_N$ where $V_N$ is a reduced basis of dimension $N=110$ that has been generated by running the classical greedy algorithm from \cite{BMPPT} over a training set $\cM_{\greedy}$ of $10^3$ snapshots. We recall that {\wnew an idealized version} is defined for $n\geq 1$ as
\begin{equation}
\label{eq:greedyAlg}
u_n \in \argmax_{u \in \cM_{\greedy}} \Vert u - P_{V_{n-1}} u\Vert,\quad V_n \coloneqq V_{n-1}\oplus \R u_n={\rm span}\{u_1,\dots,u_n\},
\end{equation}
with the convention $V_0\coloneqq \{0\}$. Figure \ref{fig:greedy} gives the decay of the error
$$
e_n^{\text{(greedy)}} = \max_{u \in \cM_{\text{greedy}}} \Vert u - P_{V_n}u \Vert
$$
across the greedy iterations.

We next estimate the truncation accuracy $\e_N$ defined in \eqref{eq:epsN}. This has been done by computing the maximum 
of the error $\Vert u - P_{V_N} u \Vert$ over the training set $\cM_{\greedy}$
supplemented by a test set $\cM_{\test}$, also of $10^3$ snapshots. We obtain the estimate
$$
\max_{u\in \cM_{\greedy} \cup \cM_{\test}} \Vert u - P_{V_N} u \Vert \leq \e_N = 5.10^{-5}.
$$
In the comparison of the three different reconstruction algorithms, we want to illustrate the impact of the number of measurements that are used. To do this, we consider the nested subspaces
$$
W_m = \Span\{ \omega_{\bx_i,\tau_i} \}_{i=1}^m \subset W
$$
for $m = 10,\,20,\,30,\,40,\,50$ so that $W_{50}=W$.

For the computation of the best affine algorithm, we generate a new training set $\wt\cM$ of $10^3$ snapshots which we project into $W+V_N$. This projected set, which we denote by $P_{W+V_N}\wt\cM$ with a slight abuse of notation, is used to compute
$$
\wt A^{(m)}_{\wca}(u) =  \wt c^{(m)} + \wt B^{(m)} P_{W_m} u,\quad m= 10,\, 20,\dots,\,50,
$$
by running the primal-dual algorithm of section \ref{sec:optim}. We have added the indices $m$ to stress that the algorithm depends on it.

%
%

\begin{figure}
    \centering
    \begin{subfigure}[b]{0.45\textwidth}
        \includegraphics[width=\textwidth]{{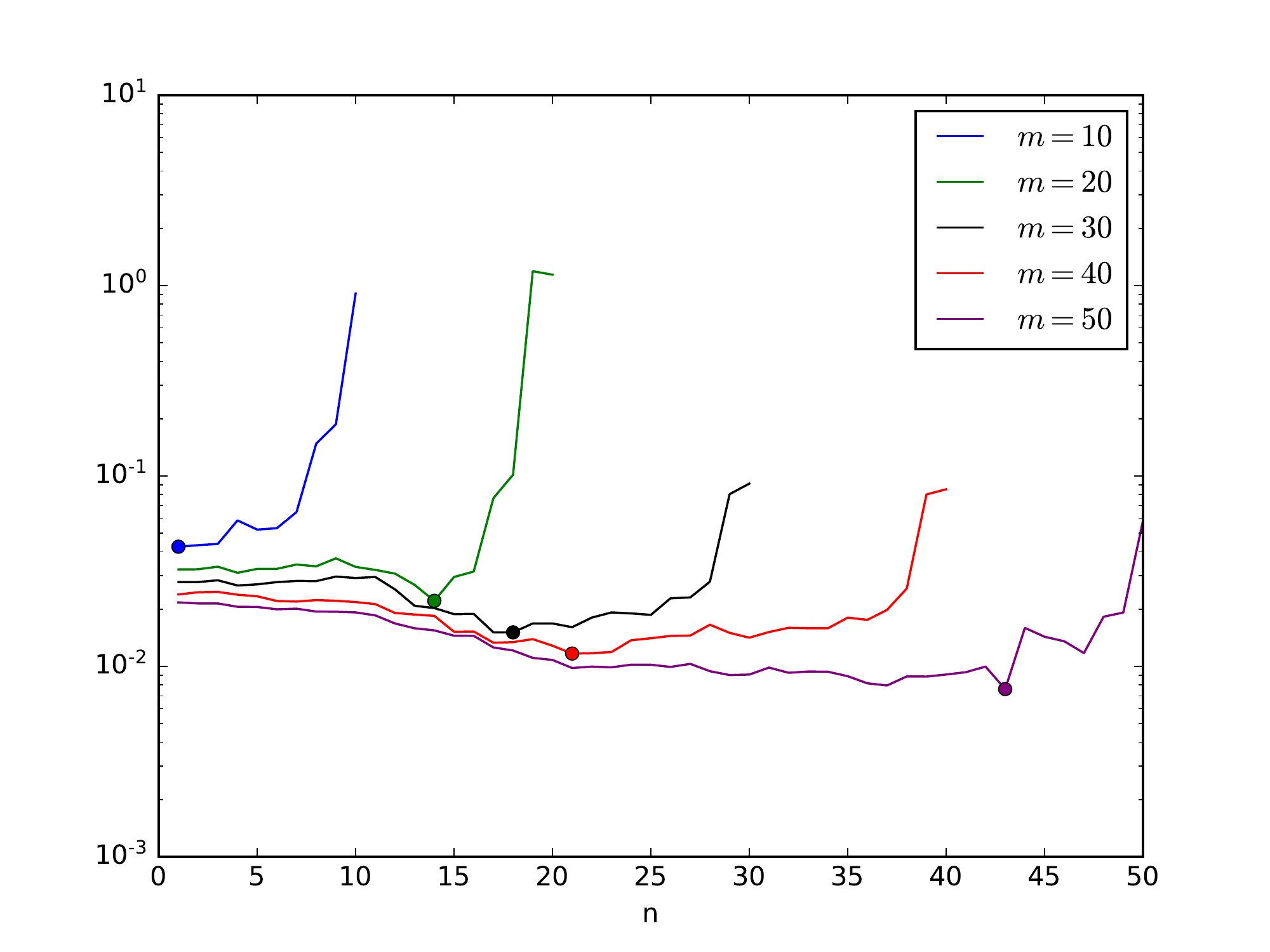}}
        \caption{$e_{\one}^{(m,n)}$: reconstruction error with $A^{(m)}_n$.}
        \label{fig:one-space-Amn}
    \end{subfigure}
    ~
    \begin{subfigure}[b]{0.45\textwidth}
        \includegraphics[width=\textwidth]{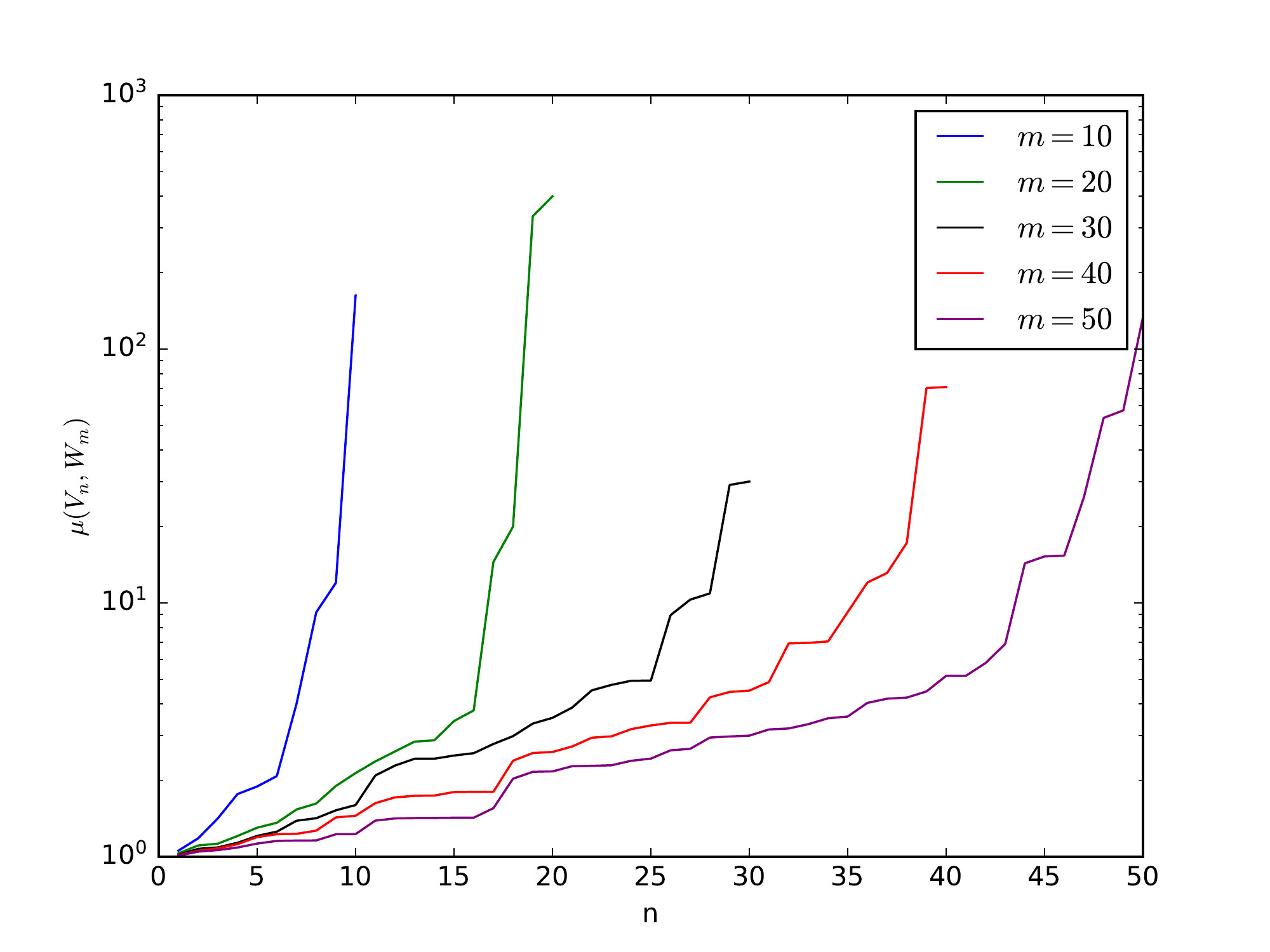}
        \caption{Stability factor $\mu(\wt V_n, W_m)$.}
        \label{fig:gull}
    \end{subfigure}    
\caption{one-space algorithm.}
\label{fig:one-space}
\end{figure}

For the comparison with the three other reconstruction algorithms, we evaluate
$$
e_{\wca}^{(m)} =\max_{u\in \cM_{\test}} || u - \wt A^{(m)}_{\wca}(P_{W_m} u) ||
,\quad m= 10,\, 20,\dots,\,50.
$$
We stress on the fact that the three sets $\cM_{\greedy},\ \wt\cM$ and $\cM_{\test}$ are {\it different}.
 We compare this value with the performance of a straightforward reconstruction with the minimal $V$-norm recovery map,
$$
e_{\mvn}^{(m)} =\max_{u\in \cM_{\test}} || u - P_{W_m} u ||
,\quad m= 10,\, 20,\dots,\,50,
$$
with the mean square approach,
$$
e_{\msa}^{(m)} =\max_{u\in \cM_{\test}} || u - \wt A^{(m)}_{\msa}(P_{W_m} u) ||
,\quad m= 10,\, 20,\dots,\,50,
$$
and with the best {\wnew one-space} affine algorithm,
$$
e_{\one}^{(m)} = \min_{1\leq n \leq m} e_{\one}^{(m,n)},
$$
where
\be
e_{\one}^{(m,n)}
=
\max_{u\in \cM_{\test}} || u - A^{(m)}_n (P_{W_m} u) ||
,\quad m= 10,\, 20,\dots,\,50.
\label{eq:err-one-mn}
\ee
Some remarks on the computation of the {\wnew one-space} algorithm are in order.
First of all, we have used the average
$$
\bar u := \frac{1}{\#  \cM_{\greedy}} \sum_{u \in \cM_{\greedy}} u
$$
as our offset. For $m\leq M$ and $n \leq m$ given, the {\wnew one-space} affine algorithm  $A^{(m)}_n$ is the one involving the spaces $W_m$ and $\wt V_n = \bar u + V_n$, where $V_n={\rm span}\{u_1,\dots,u_n\}$. Its performance is given by $e_{\one}^{(m,n)}$ in formula \eqref{eq:err-one-mn}. Figure \ref{fig:one-space-Amn} shows $e_{\one}^{(m,n)}$ as a function of $n$ and $m$. Note that, for a fixed $m$, the error $e_{\one}^{(m,n)}$ reaches a minimal value $e_{\one}^m = e_{\one}^{(m,n^*)}$ for a certain dimension $n^*=n^*(m)$ of the reduced model, given by a thick dot in the figure. This behavior is due to the trade-off between the increase of the approximation properties of $\wt V_n$ as $n$ grows and the degradation of the stability of the algorithm, given by the increase of $\mu(\wt V_n, W_m)$ with $n$. For our comparison purpose, we use $A^{(m)}_{\one} = A_{n^*(m)}^{(m)}$, that is,  
the best possible {\wnew one-space} algorithm based on the reduced basis spaces.

\begin{figure}
\centering
\includegraphics[scale=0.4]{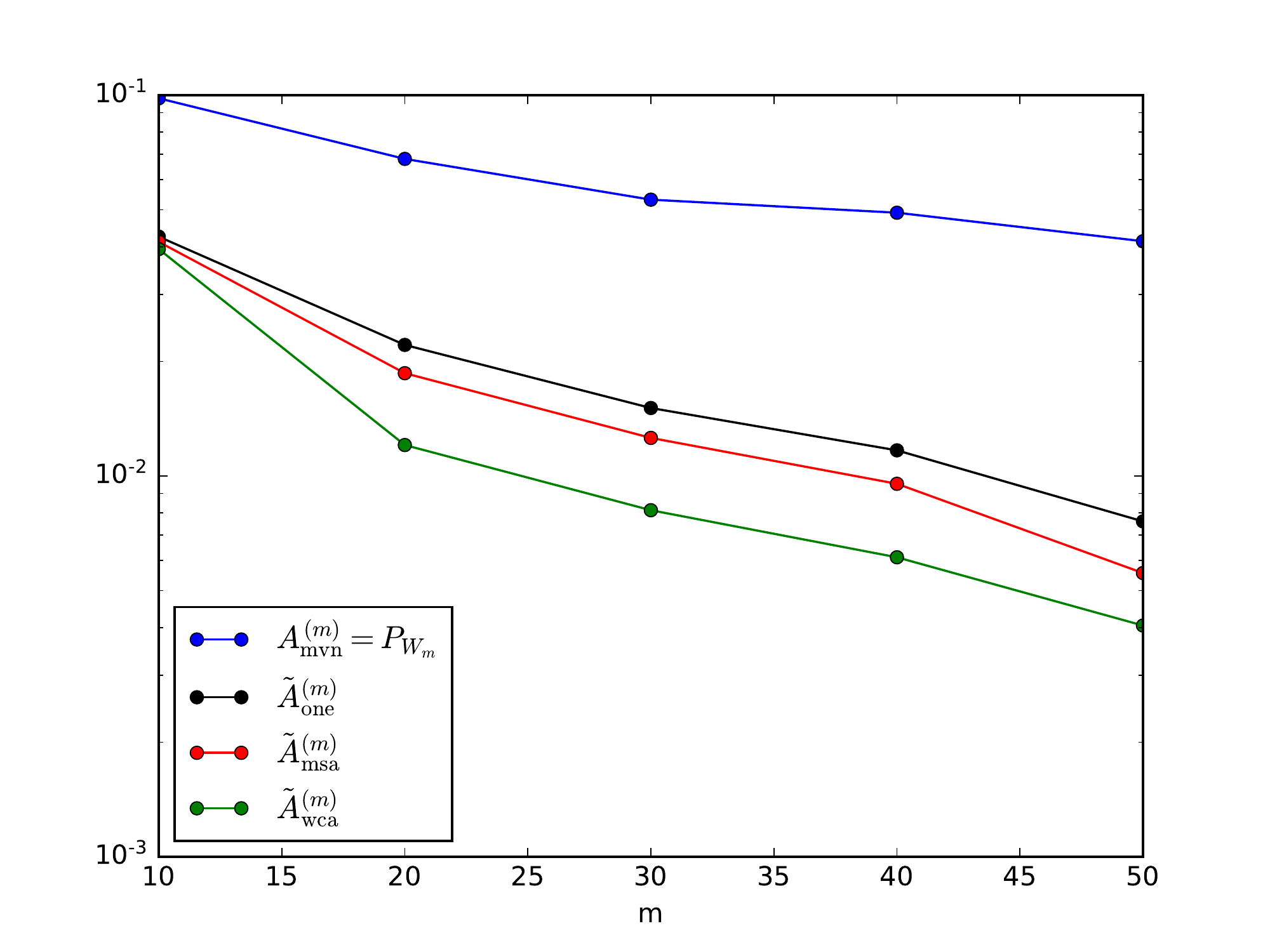}
\includegraphics[scale=0.4]{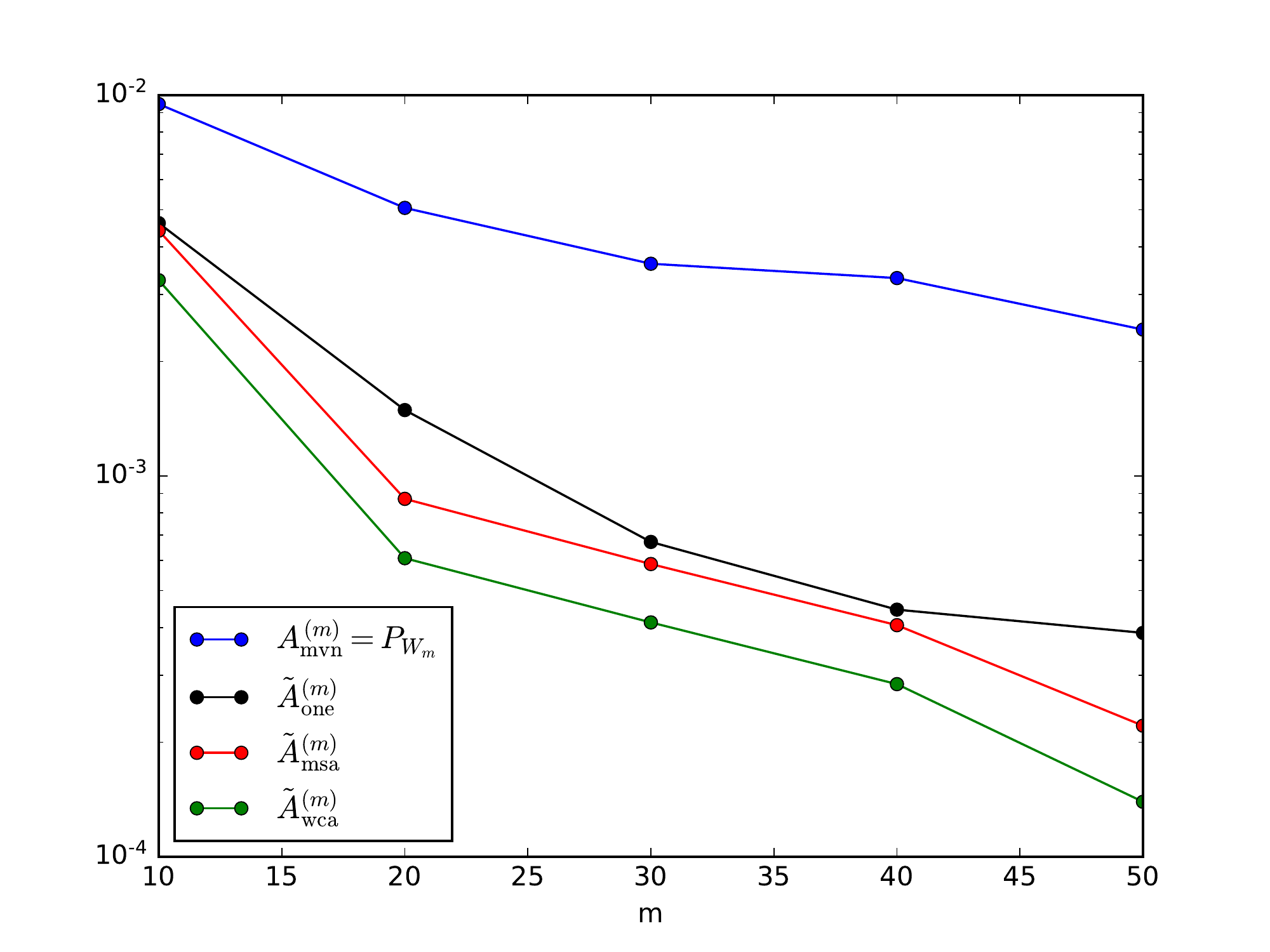}
\caption{Comparison of the reconstruction errors (left: $H^1_0(D)$ norm; right: $L^2(D)$ norm).}
\label{fig:err-summary}
\end{figure}

Figure \ref{fig:err-summary} shows the reconstruction errors $e_{\wca}^{(m)} $, $e_{\mvn}^{(m)}$, $e_{\msa}^{(m)}$ and $e_{\one}^{(m)}$ of the four different approaches for $m=10,\,20,\dots,\,50$. We also append a table with the values. We observe that a straightforward reconstruction with the minimal $V$-norm algorithm performs poorly in terms of approximation error and its quality improves only very mildly as we increase the number $m$ of measurements. This justifies considering our three other, more sophisticated, reconstruction algorithms. In this respect, the results confirm first of all that $\wt A^{(m)}_{\wca}$ is the best reconstruction algorithm. The mean square approach appears to be slightly superior to the {\wnew one-space} algorithm but still worse than the best affine algorithm. Note  that the accuracy improvement between the best affine algorithm and the {\wnew one-space} and mean square algorithms is of about a half order of magnitude for each $m$.

\begin{figure}
\centering
\includegraphics[width=0.45\textwidth]{{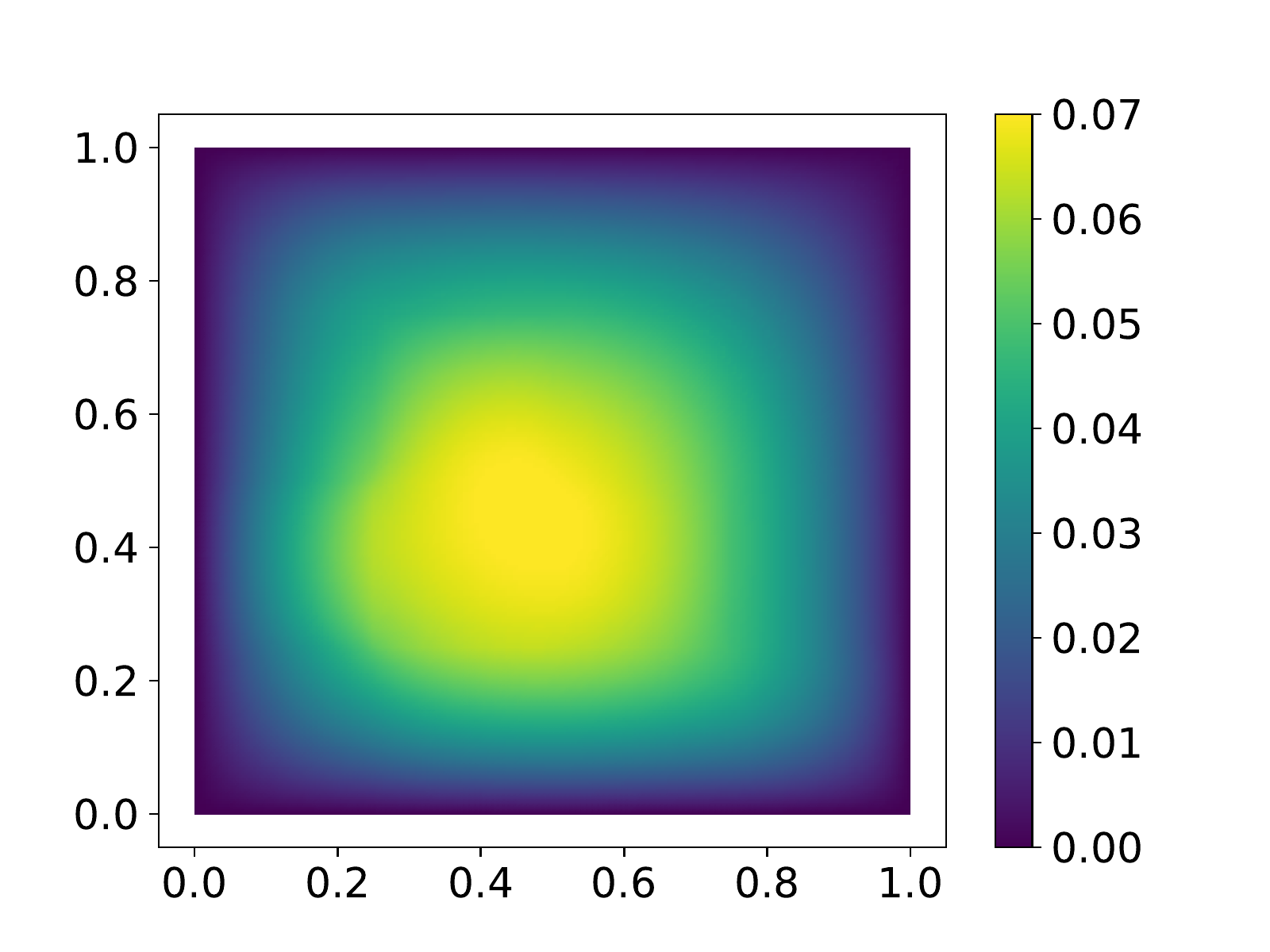}}
\includegraphics[width=0.45\textwidth]{{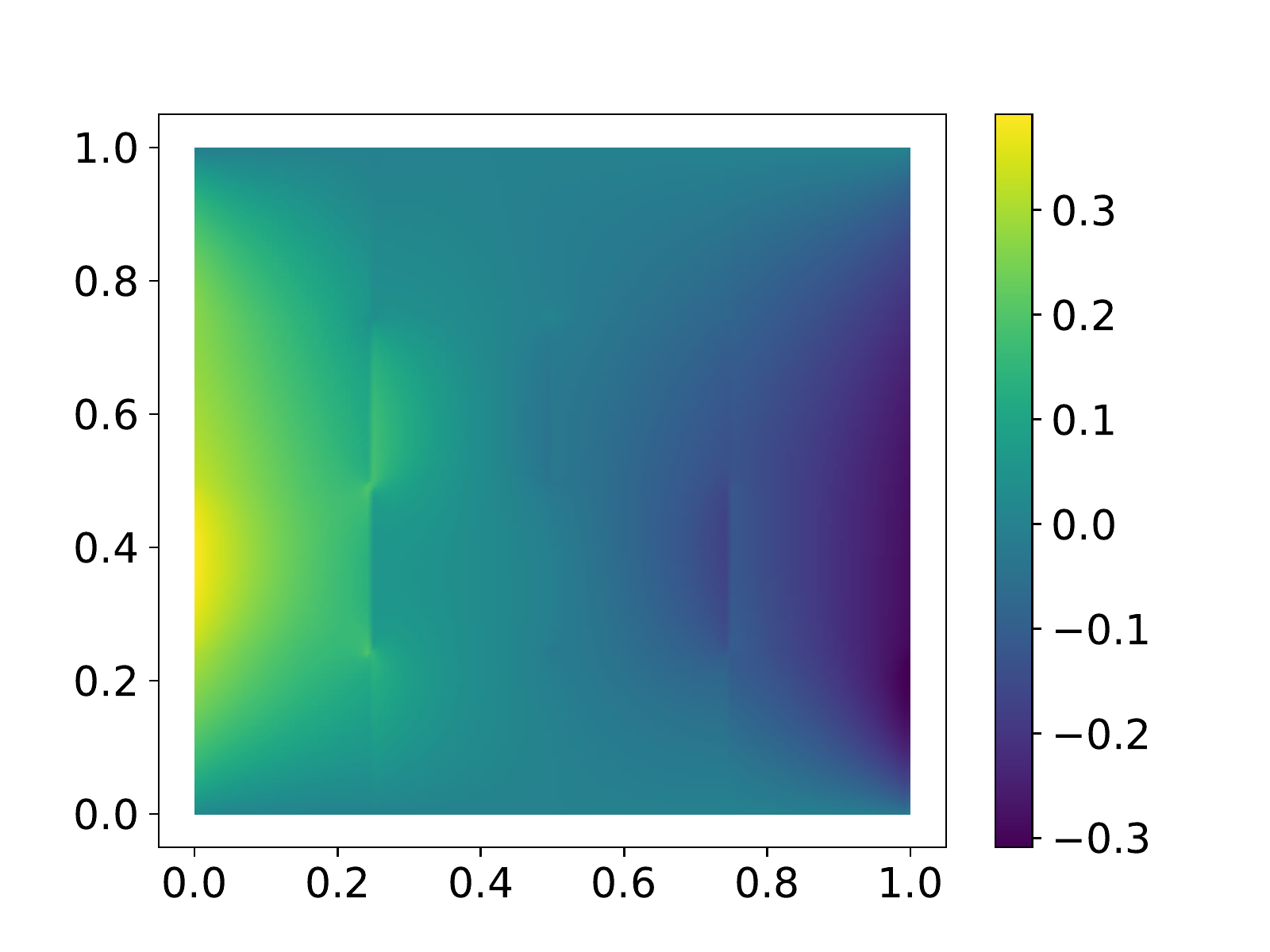}}
\caption{Function $u$ (left) and $\partial u/\partial x$ (right). The reconstruction of this function is given in Figures \ref{fig:img-reconstruction-m-20} and \ref{fig:img-reconstruction-m-40}.}
\label{fig:snap-exact}
\end{figure}

Last but not least, we give some illustrations on the reconstruction algorithms applied to a particular snapshot function $u$ from the test set $\cM_{\test}$. The target function is given in Figure \ref{fig:snap-exact} and Figures \ref{fig:img-reconstruction-m-20} and \ref{fig:img-reconstruction-m-40} show the resulting reconstructions of $u$ from $P_{W_m}u$ with our four different algorithms and for $m=20$ and $40$. Visually, the reconstructed functions look very similar. However, the difference in quality can be better appreciated in the plots of the spatial errors $| u(\bx) - A^{(m)}(u)(\bx) |$ as well as in the derivatives and their corresponding spatial errors.

\newcommand{\makefunplots}[1]{%
\begin{figure}
    \centering
    \begin{subfigure}[b]{\textwidth}
    \centering
        \includegraphics[width=0.24\textwidth]{snap-meas-m-#1.pdf}
         \includegraphics[width=0.24\textwidth]{snap-meas-m-#1-error.pdf}
         \includegraphics[width=0.24\textwidth]{snap-meas-m-#1-dudx.pdf}
         \includegraphics[width=0.24\textwidth]{snap-meas-m-#1-error-dudx.pdf}
        \caption{Minimal $V$-norm: $P_{W_{#1}}(u)$.}
    \end{subfigure}
    \\
      \begin{subfigure}[b]{\textwidth}
      \centering
        \includegraphics[width=0.24\textwidth]{snap-pbdw-aff-m-#1.pdf}
        \includegraphics[width=0.24\textwidth]{snap-pbdw-aff-m-#1-error.pdf}
        \includegraphics[width=0.24\textwidth]{snap-pbdw-aff-m-#1-dudx.pdf}
         \includegraphics[width=0.24\textwidth]{snap-pbdw-aff-m-#1-error-dudx.pdf}
        \caption{one-space affine: $A^{(#1)}_{\one}\(P_{W_{#1}}(u)\)$}
    \end{subfigure}
    \\
    \begin{subfigure}[b]{\textwidth}
    \centering
     \includegraphics[width=0.24\textwidth]{snap-msa-m-#1.pdf}
     \includegraphics[width=0.24\textwidth]{snap-msa-m-#1-error.pdf}
     \includegraphics[width=0.24\textwidth]{snap-msa-m-#1-dudx.pdf}
         \includegraphics[width=0.24\textwidth]{snap-msa-m-#1-error-dudx.pdf}
      \caption{Mean Square Algorithm: $\wt A^{(#1)}_{\msa}\(P_{W_{#1}}(u)\)$}
    \end{subfigure}  
    \\
     \begin{subfigure}[b]{\textwidth}
     \centering
     \includegraphics[width=0.24\textwidth]{snap-best-aff-m-#1.pdf}
     \includegraphics[width=0.24\textwidth]{snap-best-aff-m-#1-error.pdf}
     \includegraphics[width=0.24\textwidth]{snap-best-aff-m-#1-dudx.pdf}
         \includegraphics[width=0.24\textwidth]{snap-best-aff-m-#1-error-dudx.pdf}
      \caption{Best affine: $\wt A^{(#1)}_{\wca}\(P_{W_{#1}}(u)\)$}
    \end{subfigure}    
\caption{Reconstruction of the given function $u$ $(m=#1)$. For each reconstruction strategy: (i) the two first figures are $A^{(m)}(u)(\bx)$ and the spatial errors $| u(\bx) - A^{(m)}(u)(\bx) |$, (ii) the two last figures are $\frac{\partial A^{(m)}(u)}{\partial x}(\bx)$ and the spatial errors $| \frac{\partial u}{\partial x}(\bx) - \frac{\partial A^{(m)}(u)}{\partial x}(\bx) |$.}
\label{fig:img-reconstruction-m-#1}
\end{figure}
}

\makefunplots{20}
\makefunplots{40}

{\mnew 
Let us briefly discuss the complexity of the primal-dual algorithm.
At each iteration of the algorithm, the main bottleneck is the computation of $L^*$ (equation (3.21)). It requires to do $J$ matrix-vector products with the matrices $Q_j \in \R^{N\times m(N+1)}$ and then do a summation of the resulting vectors. The cost of these operations thus increases linearly with $J$ in terms of computational time and memory ressources. In fact, the limitation in memory was the main reason to fix $J=10^3$ and not work with a larger number of training snapshots. Let us make a quick count on the cost in terms of the number of elements to store at each iteration. The matrices $Q_j$ are sparse. For each row, there are $m+1$ nonnegative coefficients. Therefore we need to store $N (m+1)$ coefficients for each matrix, therefore a total number of $J N (m+1)$ coefficients. In our case, $N=110$ was carefully fixed to guarantee that
$$
\max_{u\in \cM_{\greedy} \cup \cM_{\test}} \Vert u - P_{V_N} u \Vert \leq \e_N = 5.10^{-5}.
$$
We have $m$ ranging between $10$ and $50$. Thus the number of nonnegative elements that we have to store for each $Q_j$ ranges between $1210$ and $5610$. Therefore, taking $J=10^3$ as in our computation, we need to handle a total number of coefficients ranging between $1.21.10^6$ and $5.61.10^6$. 
}

{\mnew 
\section{Conclusions}

In this paper, we have studied the notion of a best affine recovery map for a general state estimation
problem, that is, the map $A^*_{\rm wca}$ that minimizes the worst case error $E_{\wc}(A,\cM)$ among all affine maps.
This map is the solution to a convex optimization problem. Up to the additional perturbation induced from replacing $\cM$
by a discrete training set $\wt \cM$, it can be efficiently computed by a primal-dual optimization algorithm.
Since any affine recovery map is associated with a reduced basis $V_n$ plus an offset $\bar u$, 
the optimal affine map amounts to applying the {\wnew one-space} method from \cite{MPPY}
using an affine reduced model space $\bar u +V_n$ which is optimal for the reconstruction task. Our numerical tests confirm 
that this choice outperforms standard reduced basis spaces, which are not specifically constructed 
for the recovery problem, but rather for the approximation of $\cM$.

Our approach is readily applicable to any type of parametric PDEs, ranging from linear PDEs with affine parameter dependence to 
non-linear PDEs with non-affine parameter dependence.  We outline  its  main limitations:
\begin{itemize}
\item The first essential limitation lies in its confinement to linear or affine recovery algorithms.
Let us stress that state estimation is a linear inverse problem in the sense that the observed data $w$ is generated
from $u$ by a linear projection, optimal recovery among all possible maps, due to the
complex nonlinear geometry of the solution manifold $\cM$ that constitutes the prior. Therefore, going beyond
the results provided by our method requires the development of nonlinear recovery strategies.
One possible approach, currently under investigation, is to (i) consider a collection affine reduced model spaces 
$\{\bar u_k+ V^k \, : \, k=1,\dots K\}$, each of them of dimension $n_k\leq m$, (ii)
use the observed data $w$ to properly select a particular space from this collection
and (iii) apply the affine recovery algorithm using this particular data-dependent space.
One standard way to obtain such local reduced model spaces is by splitting the parameter domain
and searching for local reduced bases or POD, as for example proposed in \cite{Ams} for forward modeling
or in \cite{LM} for state estimation. However, the
optimal affine recovery approach discussed in the present paper could also be used in order to 
improve on such constructions.
\item The developed approach implicitly assumes that the parametric PDE model is perfect although, in full generality, the true physical state may not belong to $\cM$. It is also assumed that measurements are noiseless. One way to readily extend this approach to the search of optimal affine maps that take into account model bias and measurement noise
is as follows: suppose that the model bias is of size $\delta>0$ in the sense that the real physical state $u$ belongs to the offset
$$
\cM_\delta := \{ v \in V \,:\, \exists y\in Y \text{ s.t. }  \Vert v - u(y) \Vert \leq \delta \}.
$$
Suppose further that measurements are given with some deterministic noise, that is, we are given $P_W u + \eta$ such that $\Vert \eta \Vert \leq \sigma$ for some noise level $\sigma$. Then, the optimal affine map is given by
$$
\min_{A \text{ affine}}\; \max_{\substack{u\in \cM_\delta,\\  \| \eta \|\leq \sigma}}  \|u - A(P_Wu + \eta ) \|
$$
Once again we may emulate this optimization by introducing a discrete training set. Let us stress
that such an optimization problem requires the knowledge of the size of the model bias $\delta$ 
and the noise level $\sigma$. While $\sigma$ may be known for some applications, $
\delta$ is very hard to estimate in practice.
\end{itemize}
}
{\wnew An assessment of the obtained estimation accuracy relies, however,
on the availability of computable bounds for  the distance of the reduced spaces from the solution manifold which 
may depend on the type of the PDE model.}

\end{document}